\documentclass[12pt]{article}
\DeclareFontFamily{U}{mathx}{\hyphenchar\font45}
\DeclareFontShape{U}{mathx}{m}{n}{
      <5> <6> <7> <8> <9> <10>
      <10.95> <12> <14.4> <17.28> <20.74> <24.88>
      mathx10
      }{}
\DeclareSymbolFont{mathx}{U}{mathx}{m}{n}
\DeclareFontSubstitution{U}{mathx}{m}{n}
\DeclareMathAccent{\widebar}{0}{mathx}{"73}
\usepackage[dvips]{graphicx}
\usepackage{epstopdf}
\usepackage{color}
\usepackage{amsfonts}
\usepackage{amsmath,amssymb,graphics,color,mathrsfs,dsfont,amsfonts,extpfeil}
\usepackage{indentfirst}
\usepackage{bm}
\usepackage{multirow}
\usepackage{subfigure}
\usepackage{lscape}
\usepackage{graphicx}
\usepackage{epsfig}
\usepackage{amsfonts}
\usepackage{pst-plot}
\usepackage{multirow,booktabs}
\usepackage{algorithm}
\usepackage{algorithmic}
\usepackage{hyperref}
\usepackage{setspace}
\hypersetup{
    colorlinks=true,
    linkcolor=blue,
    filecolor=blue,
    urlcolor=blue,
    citecolor=cyan,
}

\makeatletter

\newcommand{\Rmnum}[1]{\expandafter\@slowromancap\romannumeral #1@}
\makeatother

\newcommand{\no}{\nonumber}

\newcommand{\be}{\begin{equation}}
\newcommand{\ee}{\end{equation}}
\newcommand{\bea}{\begin{eqnarray}}
\newcommand{\eea}{\end{eqnarray}}
\newcommand{\beaa}{\begin{eqnarray*}}
\newcommand{\eeaa}{\end{eqnarray*}}

\newtheorem{lemma}{Lemma}
\newtheorem{theorem}{Theorem}
\newtheorem{proposition}{Proposition}

\newcommand{\bX}{\bm{X}}
\newcommand{\bG}{\bm{G}}
\newcommand{\bY}{\bm{Y}}
\newcommand{\bZ}{\bm{Z}}
\newcommand{\bs}{\bm{s}}
\newcommand{\bt}{\bm{t}}
\newcommand{\bve}{\bm{e}}
\newcommand{\bw}{\bm{w}}
\newcommand{\mH}{\mathcal{H}}
\newcommand{\vep}{\varepsilon}
\newcommand{\meta}{\bm{\eta}}
\newcommand{\bbeta}{\bm{\beta}}
\newcommand{\balpha}{\bm{\alpha}}
\newcommand{\bphi}{\bm{\phi}}
\newcommand\independent{\protect\mathpalette{\protect\independenT}{\perp}}
\def\independenT#1#2{\mathrel{\rlap{$#1#2$}\mkern2mu{#1#2}}}
\begin{document}
\baselineskip 20pt
\parskip=10pt

% \title[short text for running head]{full title}{comments for title}
\begin{center}{\bf{\LARGE A new reproducing kernel based nonlinear dimension reduction method for survival data\footnote[1]{This research is supported by NSFC. Grant No. 71873128, the
Natural Science Foundation of Anhui Province (Grant No. 1308085MA02), and the Natural Sciences and Engineering Research Council of Canada (RGPIN 2017 05720).}}}\end{center} \vspace{0.1in}

% \author[]{Full name}{footnote}
% Remark:  One \author for one author
\begin{center}

{\bf Wenquan Cui$^1${\footnote[2]{Corresponding author: wqcui@ustc.edu.cn}}\qquad Jianjun Xu$^1$ \qquad   and\qquad  Yuehua Wu$^2$}

{\em\footnotesize  $^1$Department of Statistics and Finance, University of Science and of Technology of China,Hefei, {\rm 230026}, China\\
$^2$Department of Mathematics and Statistics, York University, Toronto, Ontario {\rm M3J 1P3}, Canada}

\end{center}

\begin{spacing}{1.2}
\begin{abstract}
Based on the theories of sliced inverse regression (SIR) and reproducing kernel Hilbert space (RKHS), a new approach
RDSIR (RKHS-based Double SIR) to nonlinear dimension reduction for survival data is proposed and
discussed. An isometrically isomorphism is constructed based on RKHS property, then the nonlinear function in the RKHS can be represented by the inner product of two elements which reside in the isomorphic feature space. Due to the censorship of survival data, double slicing is used to estimate weight function or conditional survival function to adjust for the censoring bias.
The sufficient dimension reduction (SDR) subspace is estimated
 by a generalized eigen-decomposition problem. Our method is computationally efficient with fast calculation speed and small computational burden
 The asymptotic property and the convergence rate of the estimator are also discussed based on the perturbation theory.
 Finally, we illustrate the performance of RDSIR on simulated and real data to confirm that RDSIR is comparable with linear SDR method.  The most important is that RDSIR can also extract nonlinearity in survival data effectively.
 %\vskip0.2cm
 %\noindent\textbf{Keywords}:\ RKHS; sliced inverse regression; central subspace; survival data; perturbation theory\vspace{-3mm}
\end{abstract}
\end{spacing}
\vspace{-1mm}
%  Keyword is required.

%\keywords{RKHS, sliced inverse regression,  sufficient dimension reduction subspace, survival data, operator, perturbation theory}

%  \subjclass is required.
% \MSC{\bXX\bXX\bX, \bXX\bXX\bX}
\section{Introduction}
Sparse high-dimensional data are encountered in a wide range of areas including biology, genomics, health sciences, astronomy, economics and machine learning. \text{var}iable selection and sufficient dimension reduction (SDR) are two commonly used methods in modeling such data. However, these two methods are based on different assumptions. For variable selection, researchers assume that among all the covariates, only a few are truly related to the response.
In recent years, researchers have done much work on variable selection and gained notable achievements including  LASSO type \cite{LASSO,ALASSO,grouplasso,scaledlasso}, Elastic net \cite{LASSOnet}, Dantzig selector \cite{Dantzig2007}, SCAD \cite{scad2}, SIS \cite{SIS1,SIS2} and many others.
While in the aspect of SDR, the assumption is that the response variable may relate to all the covariates but only relates to a few linear combinations of them. The goal of SDR is to recover the space
spanned by the coefficient vectors of these linear combinations. There are three commonly used methods in the existing literature: the inverse regression methods, the non-parametric methods (e.g.\cite{MAVE}) and the semiparametric methods (e.g. \cite{MaZhu12}).
The inverse regression methods started with sliced inverse regression (SIR) proposed by Li~\cite{Li1991},
which considered the following model
\begin{equation}\label{DRmodel}
   Y=g(\bm{\beta}_1^T\bX,\cdots,\bm{\beta}_q^T\bX,\epsilon),
\end{equation}
where $\bX$ is the $p$-dimensional regressor vector, $Y$ is the response variable, $g$ is an arbitrary function on $\mathbb{R}^{q+1}$, random error $\epsilon$ and $\bX$ are independent, $\bm{\beta}_1,\cdots,\bm{\beta}_q$ are $p$-dimensional non-random vectors which are so called SDR directions which capture all we need to know about $Y$. The space generated by SDR directions is SDR space. Obviously, $\{\bm{\beta}_1,\cdots,\bm{\beta}_q\}=I_p$ is a trivial case.
Denote $\mathcal{B}=\mathrm{span}\{\bm{\beta}_1,\cdots,\bm{\beta}_q\}$, the central subspace \cite{COOK1998} with the smallest $q\ (q\ll p)$, which is the intersection of all~SDR subspaces. The goal of SDR under model~(\ref{DRmodel}) is to find~$\mathcal{B}$ without specifying the unknown function $g$.
{For identifying the central space, SIR uses the first several significant eigenvectors of the matrix $\mbox{Cov}\{E(\bX)|Y\}$
to recover $\mathcal{B}$. However, SIR fails to identify some symmetric patterns, sliced average variance estimation (SAVE) \cite{COOK1991} was developed to deal with this problem, but SAVE is less efficient than SIR \cite{COOK1999}.
Some other inverse regression methods are directional regression (DR) \cite{directionalregression}, kernel inverse regression \cite{kernelIR}, and so on. All these methods are based on certain conditions including linear design condition \cite{Li1991} and constant variance condition \cite{COOK1991}. As for non-parametric methods that are conceptually more intuitive while computationally more complicated, the most original ones are the
minimum average variance estimation (MAVE) \cite{MAVE} and density based MAVE (dMAVE) \cite{dMAVE}. These methods do not require the linear design condition and constant variance condition that are critical for inverse regression methods.
Ma and Zhu \cite{MaZhu12} casted the SDR problem in a semiparametric estimation framework, in which the common conditions of linearity and constant variance on the covariates may be removed but at the cost of performing nonparametric regression.
There are many other methods available in literatures. See~\cite{MaZhu13} among others for more details.

Until now, we only discuss that the reduced predictors take the linear form $\bbeta^T\bX$, which may not work if the response does not relate to a function of a finite set of $\{\bbeta^T\bX, \ \bbeta\in \mathbb{R}^p\}$. In contrast, the reduced predictors may take  the nonlinear form $u(\bX)$, where $u$ is a function in a Hilbert space. Thus, linear SDR methods were
generalized to the nonlinear ones, e.g., generalized SIR (GSIR) and generalized SAVE (GSAVE) \cite{GSIRGSAVE},
kernel SIR \cite{Wu2008}, kernel dimension reduction (KDR)\cite{Fukumizu09} and \cite{Wuqiang2013,Yeh09} and so on. But as far as we know, the existing such nonlinear methods are not suitable for survival data.}

Survival analysis is a branch of statistics that is used to analyze data in which the time until the event is of interest. Survival data are often censored, and hence it is always a challenge to model such data. Thus, the research on SDR for survival data is not very well developed. DSIR (Double SIR) proposed by Li \textit{et al.} \cite{Li1999} is an extension of SIR for survival data. DSIR assumes that the survival time~$T$ and the regressor vector $\bX$ can be modeled by the model~(\ref{DRmodel}). Thus, when the linear design condition~\cite{Li1991} holds, $E(\bX|T)-E(\bX)\in \mathrm{span}\{\Sigma\bm{\beta}_j, j=1,\cdots,q\}$,  where~$\Sigma=\text{Cov}(\bX)$ is the covariance matrix of~$\bX$. For characterizing the information of censoring, DSIR assumes that the censoring time~$C$ can be modeled by a model like~(\ref{DRmodel}). However, DSIR only finds linear SDR directions. Therefore, we need new nonlinear dimension reduction approaches for survival data in order to extract important nonlinear {components}. The literature on SDR for survival data in recent decades includes  \cite{LiLi04,MayaStephan14,Xuerong10,XIA2011,lu2011} among others. However, most of them discussed only linear SDR directions. To our knowledge, nonlinear SDR for survival data is a completely new field.

 In this paper, we go one step further to propose a method for nonlinear SDR for survival analysis. Our method is motivated by the DSIR and the kernel trick that is an ingenious technique based on the reproducing kernel Hilbert space (RKHS) theory. See~\cite{RKHSfirst,Berlinet2011,Bernhard02,Wahba1990} for more details about RKHS. Specifically, we allow the nonlinear functions to reside in a RKHS, an isometrically isomorphism is constructed based on RKHS property, then $u$ above can be
represented by the inner product of two elements which reside in the isomorphic feature space. The linear directions found in the isomorphic feature space are corresponding to the nonlinear directions in the original input space.
Zhong \textit{et al.} \cite{RSIR1991} proposed the regularized SIR (RSIR) method which is efficient and stable in computation for data with high dimensionality and high collinearity. Wu \textit{et al.} \cite{Wuqiang2013} discussed some asymptotic properties of RSIR. As in this paper, we also adopt this regularization method to improve the efficiency and stability of the computation {of the proposed method}. Cui and Wu \cite{Cui15} gave an overview of the approach to a
nonlinear SDR for censored survival data and presented some preliminary theoretical and experimental results on the method which we abbreviate as RDSIR (RKHS-based Double SIR). In this paper, we give a systematic illustration on this approach and prove several theoretical results about RDSIR rigorously, which include the consistency as well as the convergence rate of nonlinear SDR directions if the regularization parameter satisfies some conditions.

The rest of the paper is organized as follows. In section~2, we introduce the nonlinear dimension reduction model for modeling survival data and then discuss the generalized eigen problem based on RKHS theory. In section 3, we investigate the asymptotic properties of RSDR. In section 4 we do Monte Carlo simulations on~RDSIR and compare it with DSIR. At last, the performance of RDSIR is illustrated on real data.

Throughout this paper, $\mathds{Z}=\left\{1,2,\ldots,\right\}$ and $\mathds{R}=(-\infty,\infty)$. The operation ``$\otimes$'' is defined as follows: for any $\balpha,\  \bbeta\in\mathcal{W}$,\ $\balpha\otimes \bbeta$ means for any $\meta\in\mathcal{W}$ such that
$(\balpha\otimes \bbeta)\meta=\langle \bbeta,\meta\rangle \balpha$. Let $\mathcal{W}_1$ and $\mathcal{W}_2$ be two subspaces of $\mathcal{W}$ such that for any $\bbeta_1\in\mathcal{W}_1$ and any $\bbeta_2\in\mathcal{W}_2$, $\bbeta_1$ is orthogonal to $\bbeta_2$, then $\mathcal{W}_1$ is said to be orthogonal to $\mathcal{W}_2$. If in addition, any $\bbeta\in\mathcal{W}$ can be written as a sum of a $\bbeta_1\in \mathcal{W}_1$ and a $\bbeta_2\in \mathcal{W}_2$, then we write $\mathcal{W}=\mathcal{W}_1\oplus\mathcal{W}_2$. Let $\mathcal{S}$ be a subspace in a Hilbert space $\mathcal{W}$. $\mathbf{P}_{\mathcal{S}}$ denotes the orthogonal projection operator on $\mathcal{S}$ such that for any $\bbeta\in\mathcal{W}$, $\mathbf{P}_{\mathcal{S}}\bbeta\in \mathcal{S}$ and is orthogonal to $\bbeta-\mathbf{P}_{\mathcal{S}}\bbeta\in\mathcal{W}$. Let $\bX$, $\bY$ and $\bZ$ be random vectors. The notation $\bX\independent\bY$ means that $\bX$ and $\bY$ are independent, while the notation $\bX\independent\bY \mid \bZ$ means that $\bX$ and $\bY$ are conditionally independent given $\bZ$. $I_n$ denotes the $n\times n$ identity matrix, and $\bm{1}_n$ stands for the $n$ dimensional vector $(1,\ldots,1)^T$. Let $A$ be an $n\times m$ matrix with $m\leq n$. Then $P_A=A(A^TA)^{-1}A^T$ is an $n\times n$ projection matrix onto the linear space spanned by columns of $A$. Let $B$ be a set. $I(B)$ denotes the indicator function of the set $B$. If $\Omega$ is a set, $\#\{\omega\in\Omega\}$ denotes the number of elements in $\Omega$.

\section{Nonlinear sufficient dimension reduction via reproducing kernel}

\subsection{The model setup}

Let $T$, $C$ and $\bm{\bX}=(X_1,\ldots,X_p)^T\in\mathcal{\bX}\subset\mathds{R}^p$ which are random, denote respectively the true (unobservable) lifetime, the censoring time and $p$-dimensional covariates. Put $\widetilde{T}=T\wedge C$ and $\Delta=I(T\leq C)$. The commonly made assumption for the censoring mechanism is given below:

\textbf{Condition 1}: $T\independent C\mid \bX$.

Let $\mathcal{H}_R$ be a RKHS with reproducing kernel~$\mathscr{R}(\cdot,\cdot)$ whose spectral decomposition is given by
\begin{equation*}
  \mathscr{R}(\bm{s},\bm{t})=\sum_{j=1}^\infty a_j\phi_j(\bs)\phi_j(\bt), \quad\bs,\bt\in \mathcal{\bX},
\end{equation*}
where $\{a_j, \ j\in \mathds{Z}\}$ is a sequence of nonnegative, non-increasing eigenvalues, $\{\phi_j(\cdot)\ j\in \mathds{Z}\}$ is a sequence of corresponding eigenfunctions. We assume that the functional dependence of $T$ on $\bX$ is given by
\begin{equation}\label{nonModle}
  T=g(u_1(\bX),...,u_q(\bX),\vep),
\end{equation}
where $g$ is a positive function, $u_j(\cdot)\in \mathcal{H}_R$, $1\leq j\leq q$, are linearly independent in $\mathcal{H}_R$, $\vep$ is the random error with zero mean, and $\vep$ is independent with $\bX\in\mathcal{\bX}$.

Let $l_2$ be the space of all sequences $x=(x_1,x_2,\cdots)$ of real numbers satisfying the condition
$\sum_{i=1}^{\infty}x_i^2<\infty.$
Define the map $\bphi:\ \mathcal{\bX}\mapsto l_2$ by
\[ \bphi(\bs)=(\sqrt{a_1}\phi_1(\bs),\sqrt{a_2}\phi_2(\bs),\cdots)^T\in l_2, \quad \bs\in\mathcal{\bX}\]
and the $\ell_2$ inner product between $\bphi(\bs)$ and $\bphi(\bt)$ is
\begin{eqnarray}
  &&\langle\bphi(\bs),\bphi(\bt) \rangle=\mathscr{R}(\bs,\bt),\ \ \bs,\bt\in\mathcal{\bX},\label{b2.2}
\end{eqnarray}
It can be shown that the space $$\mathcal{H}=\overline{{\mathrm{span}}\{\bphi(\bs_1),\bphi(\bs_2),\cdots,\bphi(\bs_m),\ \ \bs_j\in\mathcal{\bX},\ j=1,2,\ldots,m, \ \ m\in\mathds{Z}\}}$$ with the inner product (\ref{b2.2}) is a Hilbert space, where $\overline{A}$ denotes the closure of the set $A$.

Let the map $\mathcal{M}:\mathcal{H}_R\rightarrow\mathcal{H}$ satisfy $\mathcal{M}\mathscr{R}(\bs,\cdot)=\bphi(\bs)$. By \eqref{b2.2}, it follows that
\begin{eqnarray*}\label{mapping01}
 \langle \mathcal{M}\mathscr{R}(\bs,\cdot),\mathcal{M}\mathscr{R}(\bt,\cdot)\rangle=\langle\bphi(\bs),\bphi(\bt)\rangle=\mathscr{R}(\bs,\bt)=\langle \mathscr{R}(\bs,\cdot),\mathscr{R}(\bt,\cdot)\rangle_{\mathcal{H}_R},
\end{eqnarray*}
which implies that $\mathcal{M}$ is an isomorphic map from $\mathcal{H}_R$ to $\mathcal{H}$.
Hence, for each $\meta\in\mathcal{H}$, there exists a unique $u(\cdot)\in \mathcal{H}_R$ such that $\mathcal{M}u(\cdot)=\meta$. By the reproducing property of RKHS, we have
$$u(\bX)=\langle u(\cdot),\mathscr{R}(\bX,\cdot)\rangle_{\mathcal{H}_R}=\langle\bm{\meta},\bphi(\bX)\rangle.$$
Therefore, for $u_1(\cdot),\ldots, u_q(\cdot)$ in \eqref{nonModle}, there exist respectively $\meta_{1},\ldots,\meta_{q}$ such that $\mathcal{M}u_j(\cdot)=\meta_{j}$,  $1\leq j\leq q$. Then Model~(\ref{nonModle}) can be written as
\begin{equation}\label{Kernelmodel}
  T=g(\langle\bm{\eta}_{1},\bphi(\bX)\rangle,\cdots,\langle\bm{\eta}_{q},\bphi(\bX)\rangle,\varepsilon).
\end{equation}
It is easy to see that $\meta_{1},\ldots,\meta_{q}$ are the sufficient dimension reduction (SDR) directions in~$\mathcal{H}$. The linear space generated by $\meta_{1},\ldots,\meta_{q}$ is a dimension reduction space denoted by $\mathcal{B}_T\subset\mathcal{H}$. Since $u_1(\cdot),\cdots,u_q(\cdot)$ are linearly independent in $\mathcal{H}_R$, $\{\bm{\eta}_1,\cdots,\bm{\eta}_q\}$ forms a basis for $\mathcal{B}_T$, and $\mathcal{B}_T$ is hence the central dimension reduction subspace, the smallest dimension-reduction subspace. Let $\{\bbeta_1,\ldots,\bbeta_q\}$ be a basis for $\mathcal{B}_T$.
In view of \eqref{Kernelmodel}, it follows that there exists a positive function $g_1$ such that
\begin{equation*}\label{Kernelmodel2}
  T=g(\langle\bm{\eta}_{1},\bphi(\bX)\rangle,\cdots,\langle\bm{\eta}_{q},\bphi(\bX)\rangle,\varepsilon)
  =g_1(\langle\bm{\beta}_{1},\bphi(\bX)\rangle,\cdots,\langle\bm{\beta}_{q},\bphi(\bX)\rangle,\varepsilon),
\end{equation*}
which implies that
\begin{equation}\label{independentkernel}
  \bX \independent T\mid \langle \bbeta_1,\bphi(\bX)\rangle,\langle \bbeta_2,\bphi(\bX)\rangle\cdots,\langle \bbeta_q,\bphi(\bX)\rangle.
\end{equation}
We can see that the dimension reduction directions are not identifiable and
we need only to estimate the space spanned by $\{\meta_{1},\ldots,\meta_{q}\}$.

Let $\bZ$ be a random element taking values in $\mathcal{H}$. If $E\|\bZ\|^2<\infty$, the expectation of $\bZ$ denoted as $E(\bZ)\in\mathcal{H}$, is defined such that $\langle \meta,E(\bZ)\rangle=E\langle \meta,\bZ\rangle~\mbox{for~all}\ \meta\in\mathcal{H}$, and the covariance
of~$\bZ$ denoted as $\mbox{Cov}(\bZ)$, is given by $E[(\bZ-E(\bZ))\otimes (\bZ-E(\bZ))]$.

By Lemma~\ref{lemma3.1} given in Section 3, a basis $\{\bbeta_1,\ldots,\bbeta_q\}$ of $\mathcal{B}_T$ can be obtained by solving the following generalized eigenvalue-eigenvector problem:
\begin{equation}\label{eq7}
  \Gamma\bm{\beta} = \lambda\Sigma\bm{\beta},
\end{equation}
where~$\Sigma$ is the covariance operator of $\bphi(\bX)$ and $\Gamma$  is the covariance operator of $\bphi(\bX)$ conditional on~$T$, i.e.,~$\Sigma=\mathrm{Cov}(\bphi(\bX))$, $\Gamma=\mathrm{Cov}(E(\bphi(\bX)\mid T))$. From \cite{Baker1981}, $\Sigma$ is a compact operator and the
spectral decomposition of $\Sigma$ is
\begin{equation*}\label{eqqqqq}
  \Sigma=\sum_{j=1}^{\infty}\nu_j\psi_j\otimes\psi_j,
\end{equation*}
where $\{\nu_1,\nu_2,\cdots\}$ are eigenvalues and $\{\psi_1,\psi_2,\cdots\}$ are the corresponding
eigenfunctions. Under the linear design condition about the nonlinear dimension reduction as in Lemma 1,
$\Gamma$ is also a compact operator with rank $d_{\Gamma}\leq q$. The spectral decomposition of $\Gamma$ is
\begin{equation*}\label{eqqqqqq}
  \Gamma=\sum_{j=1}^{d_{\Gamma}}\zeta_j\varphi_j\otimes\varphi_j,
\end{equation*}
$\{\zeta_1,\zeta_2,\cdots\}$ and $\{\varphi_1,\varphi_2,\cdots\}$ are corresponding sequence of eigenvalues and
eigenfunctions.

\subsection{Estimating the central subspace}\label{ecs}

Partition $[0,\infty)$ by $0=t_1<t_2<...<t_L<t_{L+1}=\infty$.
Denote $\bm{\mu}=E[\bphi(\bX)]$, $\bm{\mu}_\ell=E[\bphi(\bX)|T\in D_\ell]$\ and\  $p_\ell=P(T\in D_\ell)$, where~$D_\ell=[t_\ell,t_{\ell+1})$, $\ell=1,\ldots, L$.  Throughout the rest of this paper, we assume that $p_\ell>0$, for $\ell=1,\ldots, L$.

Let $\{\widetilde{T}_i,\Delta_i,\bX_i, \ i=1,\ldots,n\}$ be the observed data of $(\widetilde{T},\Delta,\bX)$ so that $\{\widetilde{T}_i,\Delta_i,\bX_i\}, \ i=1,\ldots,n$, are independently and identically distributed (i.i.d.). In light of \eqref{eq7}, we can obtain an estimation of a basis of $\mathcal{B}_T$ by solving the following generalized eigenvalue-eigenvector problem
\begin{eqnarray}
&&  \widehat{\Gamma}{\bm{\beta}} = {\lambda}\widehat{\Sigma}{\bm{\beta}},\label{genalesti}
\end{eqnarray}
where $\widehat{\Gamma}$ and $\widehat{\Sigma}$ are respectively estimators of ${\Gamma}$ and ${\Sigma}$.
Thus, the remaining problem is to find estimations of ${\Gamma}$ and ${\Sigma}$.

{By the invariance property coined by Cook (\cite{COOK1998}), we can assume without loss of generality that $\sum_{i=1}^{n}\bm{\phi}(\bm{X}_i)={\bf 0}$. Denote the centered mapped data by $\Phi=(\bphi(\bX_1),\ldots,\bphi(\bX_n))$. We may estimate ${\Gamma}$ and ${\Sigma}$ by
\begin{eqnarray}
&&\widehat{\Sigma}= \frac{1}{n}\sum_{i=1}^n\bphi(\bX_i)\otimes\bphi(\bX_i)= \frac{1}{n}\Phi\Phi^T, \label{nnn}\\
&& \widehat{\Gamma}= \sum_{\ell=1}^L\widehat{p}_\ell\widehat{\bm{\mu}}_\ell\otimes\widehat{\bm{\mu}}_\ell= \frac{1}{n}\Phi{{Q}}\Phi^T,\label{nnnn}
\end{eqnarray}
where~${{Q}}=(1/n)\widehat{W}\widehat{P}^{-1}\widehat{W}^T$, ~$\widehat{W}=({\widehat{w}}_{il})_{1\le i\le n,1\le l\le L}$, ~$\widehat{P}=\text{diag}(\widehat{p}_1,\cdots,\widehat{p}_L)$ and $\widehat{\bm{\mu}}_\ell$ is the empirical estimator of $\bm{\mu}_\ell$. The expressions of~$\widehat{\bm{\mu}}_\ell$, $\widehat{w}_{i\ell}$ and $\widehat{p}_\ell$ will be given in Subsection \ref{douest}.

Since $\mathcal{B}_T$ is the central SDR subspace, by \eqref{nnn} and \eqref{nnnn}, the solution of \eqref{genalesti} should be in the space spanned by $\{\bphi(\bX_1),\ldots,\bphi(\bX_n)\}$, more details can be found in Lemma \ref{equasolve} given later. If $\widetilde{\bbeta}$ is such a solution, there exists an $n$-dimensional $\widetilde{\balpha}$ such that $\widetilde{\bbeta}=\Phi\widetilde{\balpha}$. Denote ${{R}}=\Phi\Phi^T$. Thus, the infinite dimensional equation (\ref{genalesti}) can be converted to the following finite dimensional one
\begin{equation}\label{equasolve1}
  {{R}}{{Q}}{{R}}{\bm{\alpha}}={\lambda}{{R}}^2{\bm{\alpha}}.
\end{equation}
Note that ${{R}}$ is singular in general. However, the relationship between the kernel matrix of centered mapped data and that of raw mapped data is known. As discussed in $\cite{Wu2008}$, this relationship is
 \begin{equation}\label{RRRR}
  {{R}}=\left(I_n-\dfrac{1}{n}1_n1_n^T\right)\widetilde{R}\left(I_n-\dfrac{1}{n}1_n1_n^T\right),
\end{equation}
where $\widetilde{R}=\mathscr{R}(\bm{X}_i,\bm{X}_j)_{1\leq i,j\leq n}$ is the Gram matrix of $\mathscr{R}(\cdot,\cdot)$ over $\bX_1,\ldots,\bX_n$.

Similar to the methods in \cite{Li2008} and \cite{RSIR1991}, we add a regularization term to the right hand side of \eqref{equasolve1} to prevent over-fitting and numerically instability, which results in
\begin{equation}\label{eq10}
 {{R}}{{Q}}{{R}}{\bm{\alpha}}={\lambda}({{R}}^2+n^2\tau{I}_n){\bm{\alpha}},
\end{equation}
where~$\tau$ is a tuning parameter. We elaborate how to chose $\tau$ in section \ref{tuningp}.}

\subsection{Double slicing}\label{douest}

$\{\widetilde{T}_i\}$ consists of two types of observations. One type of observation is the observed lifetime while the other type of observation is only the censoring time. Since the censoring time is dependent on the lifetime, to adjust for the censoring bias, we use the double slicing here in light of \cite{Li1999}.

By the definition of $\bm{\mu}_\ell$ given in Subsection \ref{ecs}, we have
  \begin{equation}\label{eq9}
    \bm{\mu}_\ell=\frac{E\left[\bphi(\bX)I\left(T\in D_\ell\right)\right]}{P\left\{T\in D_\ell\right\}}=\frac{E\left[\bphi(\bX)I\left(T\geq t_\ell\right)\right]-E\left[\bphi(\bX)I\left(T\geq t_{\ell+1}\right)\right]}{E\left[I\left(T\geq t_\ell\right)\right]-E\left[I\left(T\geq t_{\ell+1}\right)\right]}.
  \end{equation}
Since the observations of $T$ are right truncated because of censoring, we cannot directly estimate ~$E[\bphi(\bX)I(T\geq t)]$ and~$E[I(T\geq t)]$ empirically. But by Lemma~\ref{conditionsurv}, their estimators can be constructed as follows:
\begin{eqnarray}
 \widehat{E}[\bphi(\bX)I(T\geq t)] &=&  \frac{1}{n}\sum_{i:\widetilde{T}_i\geq t}^n\bphi(\bX_i)+\frac{1}{n}\sum_{i:\widetilde{T}_i<t,\;\Delta_i=0}^n\bphi(\bX_i)\widehat{w}(\widetilde{T}_i,t,\bphi(\bX_i)),\label{equ:estiphi}\no \\
  \widehat{P}\{T\geq t\} &=& \frac{1}{n}\sum_{i=1}^nI(\widetilde{T}_i\geq t)+\frac{1}{n}\sum_{i:\widetilde{T}_i<t,\Delta_i=0}^n\widehat{w}(\widetilde{T}_i,t,\bphi(\bX_i)), \label{equ:estiP}
\end{eqnarray}
where $\widehat{w}(\cdot,\cdot,\cdot)$ is an estimator of the weight function $w(\cdot,\cdot,\cdot)$ defined by
\begin{equation}\label{ww}w(t',t,\bphi(\bX))=S_T(t|T>t',\bphi(\bX))\quad \mbox{with $ t'<t$},\end{equation}
which implies that
\begin{eqnarray}
  && \widehat{E}\left[\bphi(\bX)I\left(t_\ell\leq T <t_{\ell+1}\right)\right]=  \displaystyle\frac{1}{n}\Phi\widehat{\bw}_\ell,\label{2.11b}\\
  && \widehat{p}_\ell = \widehat{P}\{T\geq t_\ell\}-\widehat{P}\{T\geq t_{\ell+1}\}= \frac{1}{n}\textbf{1}_n^T \widehat{\bw}_\ell.\label{2.12b}
\end{eqnarray}
By (\ref{eq9})$\sim$(\ref{2.12b}), an estimator of $\bm{\mu}_\ell$ is given by
\begin{eqnarray}
  \widehat{\bm{\mu}}_\ell=\frac{\widehat{E}\left[\bphi(\bX)I\left(t_\ell\leq T <t_{\ell+1}\right)\right]}{\widehat{p}_\ell}=\displaystyle\frac{\displaystyle\frac{1}{n}\Phi\widehat{\bw}_\ell}{\widehat{p}_\ell}=\frac{\Phi\widehat{\bw}_\ell}{\textbf{1}_n^T \widehat{\bw}_\ell}.\label{2.12c}
\end{eqnarray}

In the following, we present a method for estimating $w$. To make full use of censoring information, we assume that the relationship between $C$ and~$\bX$ can be modeled by
\begin{equation}
  C=h(v_1(\bX),\ldots,v_c(\bX),\epsilon),\label{modelCensoring}
\end{equation}
where~$h$  is an arbitrary function, $v_j(\cdot), \ 1\leq j\leq c$, are linearly independent in $\mathcal{H}_R$, and the random error~$\epsilon$ and the covariates $\bX$ are independent. Similar to (\ref{Kernelmodel}), Model (\ref{modelCensoring}) is expressed as
\begin{equation}\label{Kernelmodel02}
  C=h(\langle\bm{\gamma}_1,\bphi(\bX)\rangle,\cdots,\langle\bm{\gamma}_c,\bphi(\bX)\rangle,\epsilon),
\end{equation}
where $\bm{\gamma}_j\in\mathcal{H}$, $j=1,\ldots,c$. Denote ${\mathcal{B}_C}$ as the central dimension reduction subspace spanned by $\bm{\gamma}_j\in\mathcal{H},\ j=1,\ldots,c$. Let ${\mathcal{B}_J}\subset \mathcal{H}$ be the smallest Hilbert space containing both ${\mathcal{B}_T}$ and ${\mathcal{B}_C}$. Then the dimension of ${\mathcal{B}_J}$ is at most $q+c$. By \eqref{Kernelmodel} and \eqref{Kernelmodel02}, it follows that
\[(T,C)\independent X\mid \mathbf{P}_{\mathcal{B}_J}\bphi(\bX).\]
Thus, $w(t',t,\bphi(\bX))$ given in \eqref{ww} can be rewritten as
\begin{equation}\label{omegaPb}
 w(t',t,\bphi(\bX))=w(t',t,\mathbf{P}_{\mathcal{B}_J}\bphi(\bX)).
\end{equation}
Since~$E(\bphi(\bX)|\widetilde{T},\Delta)=E(E(\bphi(\bX)|T,C)|\widetilde{T},\Delta)$, under the condition of Lemma~\ref{lemma3.1}, it can be shown in the same way as in Lemma~\ref{lemma3.1} that
$$E(\bphi(\bX)|\widetilde{T},\Delta)-E(\bphi(\bX))\in \mathrm{span}\{\Sigma\bm{\beta}_1,\ldots,\Sigma\bm{\beta}_q,\Sigma\bm{\gamma}_1,\ldots,\Sigma\bm{\gamma}_c\}.$$
Thus a basis of $\mathcal{B}_J$ can be found by solving the following generalized eigenvalue-eigenvector problem
\begin{equation}\label{eq15}
  \Gamma_{J}\bm{\beta} = \lambda_J\Sigma\bm{\beta},
\end{equation}
where~$\Gamma_{J}=\mathrm{Cov}(\bm{\psi}_J(\widetilde{T},\Delta))$ with $\bm{\psi}_J(\widetilde{T},0)= E(\bphi(\bX)|\widetilde{T},\Delta=0)$ and $\bm{\psi}_J(\widetilde{T},1)= E(\bphi(\bX)|\widetilde{T},\Delta=1)$.

Similar to \cite{Li1999}, to obtain the sample version of \eqref{eq15}, we partition the interval $[0,\;\infty)$ for the censored and uncensored lifetime as follows:
$$0=t^{[\iota]}_{1}<t^{[\iota]}_{2}<\ldots<t^{[\iota]}_{L_{\iota}}<t^{[\iota]}_{L_{\iota}+1}=\infty,\quad \iota=0,1,$$
where $t^{[\iota]}_{1}<\ldots<t^{[\iota]}_{L_{\iota}}$ are respectively partitions of the interval $(0,\;\infty)$, $L_{\iota}$ denotes the number of partitions, and $\iota=0,1$ refer to the two different partitions to be applied to the censored and uncensored lifetime. Put $D^{[\iota]}_{\ell}=[t^{[\iota]}_{\ell},t^{[\iota]}_{\ell+1})$, and $n_{\iota,\ell}=\#\{i: \Delta_i=\iota, \widetilde{T}_i\in D^{[\iota]}_{\ell}\}$. The corresponding ordered time in $D^{[\iota]}_{\ell}$ consist of $\widetilde{T}^{[\iota]}_{1,\ell},\cdots,\widetilde{T}^{[\iota]}_{n_{\iota,\ell},\ell}$. Denote ${\Phi}^{[\iota]}_{\ell}=\left(\bphi(\bX^{[\iota]}_{1,\ell}),\ldots,\bphi(\bX^{[\iota]}_{n_{\iota,\ell},\ell})\right)$ and $\bm{\mu}^{[\iota]}_{\ell}=E\{\bphi(\bX)|\Delta=\iota,\widetilde{T}\in D^{[\iota]}_{\ell}\}$. The empirical version of $\bm{\mu}^{[\iota]}_{\ell}$ is given by
\begin{eqnarray*}
   \widehat{\bm{\mu}}^{[\iota]}_{\ell}=\displaystyle\frac{1}{n\widehat{p}^{[\iota]}_{\ell}}\sum_{i=1}^n \bphi(\bX_i)I(\Delta_i=\iota,\widetilde{T}_i\in D^{[\iota]}_{\ell}),\quad \iota=0,1,
\end{eqnarray*}
where~$\widehat{p}^{[\iota]}_{\ell}={n_{\iota,\ell}}/{n}$. Thus $\Gamma_{J}$ can be estimated by
\begin{equation}
\widehat{\Gamma}_{J}=\sum_{\iota\in\{0,1\}}\sum_{\ell=1}^{L_\iota}\widehat{p}^{[\iota]}_{\ell}\widehat{\bm{\mu}}^{[\iota]}_{\ell}\otimes\widehat{\bm{\mu}}^{[\iota]}_{\ell}=\frac{1}{n}\Phi{{{Q}}}_J\Phi^T,\label{JointGm01}
\end{equation}
where~${{{Q}}}_J=BAB^T$ and
\begin{eqnarray*}
&& \widetilde{\Phi}=\left({\Phi}^{[0]}_{1},\ldots,\Phi^{[0]}_{L_0},{\Phi}^{[1]}_{1},\ldots,\Phi^{[1]}_{L_1}\right)\\
&&  A=
\mbox{diag}(P_{\textbf{1}_{n_{0,1}}},\cdots,P_{\textbf{1}_{n_{0,L_0}}},P_{\textbf{1}_{n_{1,1}}},\cdots,P_{\textbf{1}_{n_{1,L_1}}}),\\
&&   B\ \mbox{is an $n\times n$ permutation matrix such that } \widetilde{\Phi}=\Phi B.
\end{eqnarray*}
Thus, an empirical version of (\ref{eq15}) is
\begin{equation}\label{eq17b}
    \widehat{\Gamma}_{J}{\bm{\beta}} = {\lambda_J}\widehat{\Sigma}{\bm{\beta}}.
\end{equation}
Denote $\mathcal{H}_1=\mbox{span}\left\{\bphi(\bX_1),\ldots,\bphi(\bX_n)\right\}$ and $\mathcal{H}_2$ to be the orthogonal complement subspace of $\mathcal{H}_1$, i.e., $\mathcal{H}_2\perp\mathcal{H}_1$ such that $\mathcal{H}=\mathcal{H}_1\oplus\mathcal{H}_2$. In view of \eqref{nnn} and \eqref{JointGm01}, for any $\bbeta_2\in\mathcal{H}_2$, $\widehat{\Gamma}_{J}{\bm{\beta}_2} = {\lambda}\widehat{\Sigma}{\bm{\beta}_2}$=0. Since we seek a basis for $\mathcal{B}_J$ which is the smallest Hilbert space containing both $\mathcal{B}_T$ and
$\mathcal{B}_C$, so the solutions of \eqref{eq17b} should be in $\mathcal{H}_1$. Denote the set of all such solutions corresponding to the nonzero $\lambda_J$ by $\{\breve{\bm{\beta}}_{J1},\cdots,\breve{\bm{\beta}}_{Jm}\}$.
%, which is then an estimated basis of ${\mathcal{B}_J}$.
Since $\breve{\bm{\beta}}_{Jl}\in\mH_1$ for any $l$, there exists $\breve{\bm{\alpha}}_{Jl}$ such that $\breve{\bm{\beta}}_{Jl}=\Phi\breve{\bm{\alpha}}_{Jl}$  for $1\leq l\leq m$. By \eqref{nnn}, \eqref{JointGm01}, and the fact that ${{R}}=\Phi^T\Phi$, it follows that
\begin{equation}\label{mmm}{{R}}{{Q}}_J{{R}}{\bm{\alpha}}={\lambda_J}{{R}}^2{\bm{\alpha}}.
\end{equation}
As done in Subsection \ref{ecs},  we modify \eqref{mmm} by using the following regularization technique
\begin{equation}\label{rr}
{{R}}{{Q}}_J{{R}}{\bm{\alpha}}={\lambda_J}({{R}}^2+n^2s{I}_n){\bm{\alpha}},\end{equation}
where~$s$ is a tuning parameter, which can be chosen by cross-validation. Denote the first $m$ significant solutions of \eqref{rr} by $\widehat{A}_J=\left(\widehat{\balpha}_{J1},\ldots,\widehat{\balpha}_{Jm}\right)$. Thus we obtain $\{\widehat{\bbeta}_{J1},\ldots,\widehat{\bbeta}_{Jm}\}$, where $\widehat{\bbeta}_{Jl}=\Phi\widehat{\balpha}_{Jl}$ for $1\leq l\leq m$.
Define ${{\bm R}}(\bX_{\le n},\bX)=(\mathscr{R}(\bX_1,\bX),\ldots,\mathscr{R}(\bX_n,\bX))^T$. We have
$$\langle\widehat{\bbeta}_{Jl},\bphi(\bX)\rangle=\langle\Phi\widehat{\balpha}_{Jl},\bphi(\bX)\rangle=\left(\Phi\widehat{\balpha}_{Jl}\right)^T\bphi(\bX)={{\bm R}}(\bX_{\le n},\bX)^T\widehat{\bm{\alpha}}_{Jl},\quad l=1,\ldots,m.$$

Define \begin{equation*}
  \Lambda(t',t|\bphi(\bX))=E\left\{\frac{I(t'\leq\widetilde{T}<t,\Delta=1)}{S_{\widetilde{T}}(\widetilde{T}|\bphi(\bX))}\Big|\bphi(\bX)\right\},
\end{equation*}
where $S_{\widetilde{T}}(\cdot|\bphi(\bX))$ is the survival function of $\widetilde{T}$ given $\bphi(\bX)$.
By Lemma \ref{conditionsurv} and \eqref{omegaPb}, it follows that
$$w(t',t,\bphi(\bX))=w(t',t,\mathbf{P}_{\mathcal{B}_J}\bphi(\bX))=\mathrm{exp}\{\Lambda(t',t|\mathbf{P}_{\mathcal{B}_J}\bphi(\bX))\}.$$
Thus, we use kernel smoothing method to estimate the conditional cumulative hazard function $\Lambda(t',t|\bphi(\bX))$ as
\begin{eqnarray}
 &&\widehat{\Lambda}(t',t|\mathbf{P}_{\mathcal{B}_J}\bphi(\bX)) =\nonumber\\
&&\frac{n^{-1}\displaystyle\sum_{\substack{i:\ t'<\widetilde{T}_i<t,\\\Delta_i=1}}^n\left[\widehat{S}_{\widetilde{T}}(\widetilde{T}_i|\mathbf{P}_{\mathcal{B}_J}\bphi(\bX_i))^{-1}h_n^{-m}K_m\left(h_n^{-1}(\bm{R}(\bX_{\le n},\bX_i)-\bm{R}(\bX_{\le n},\bX))^T\widehat{A}_{J}\right)\right]}{\widehat{f}(\mathbf{P}_{\mathcal{B}_J}\bphi(\bX))},\label{llll}\end{eqnarray}
where
\begin{eqnarray*}
&&  \widehat{S}_{\widetilde{T}}(\widetilde{T}_i|\mathbf{P}_{\mathcal{B}_J}\bphi(\bX_i))= \frac{n^{-1}\displaystyle\sum_{j:\ \widetilde{T}_j>\widetilde{T}_i}^nh_n^{-m}K_m\left(h_n^{-1}(\bm{R}(\bX_{\le n},\bX_j)-\bm{R}(\bX_{\le n},\bX_i))^T\widehat{A}_{J}\right)}{\widehat{f}(\mathbf{P}_{\mathcal{B}_J}\bphi(\bX_i))} , \\
&& \widehat{f}(\mathbf{P}_{\mathcal{B}_J}\bphi(\bX)) = n^{-1}\sum_{j=1}^nh_n^{-m}K_m\left(h_n^{-1}(\bm{R}(\bX_{\le n},\bX_j)-\bm{R}(\bX_{\le n},\bX))^T\widehat{A}_{J}\right),
\end{eqnarray*}
\iffalse
\begin{eqnarray*}
 &&\hat{\Lambda}(t',t|\mathbf{P}_{\mathcal{B}_J}\bphi(\bX)) =\\
&&\frac{n^{-1}\sum_{i:t'<\widetilde{T}_i<t,\Delta_i=1}^n(\hat{S}_{\widetilde{T}}(\widetilde{T}_i|\mathbf{P}_{\mathcal{B}_J}\bphi(\bX_i))^{-1}h_n^{-d}K_d(h_n^{-1}\mathbf{P}_{\mathcal{B}_J}(\bphi(\bX_i)-\bphi(\bX)))}{\hat{f}(\mathbf{P}_{\mathcal{B}_J}\bphi(\bX))},\\
&&  \hat{S}_{\widetilde{T}}(\widetilde{T}_i|\mathbf{P}_{\mathcal{B}_J}\bphi(\bX_i))= \frac{n^{-1}\sum_{j:\widetilde{T}_j>\widetilde{T}_i}^nh_n^{-d}K_d(h_n^{-1}\mathbf{P}_{\mathcal{B}_J}(\bphi(\bX_j)-\bphi(\bX_i)))}{\hat{f}(\mathbf{P}_{\mathcal{B}_J}\bphi(\bX))} , \\
&& \hat{f}(\mathbf{P}_{\mathcal{B}_J}\bphi(\bX)) = n^{-1}\sum_{j=1}^nh_n^{-d}K_d(h_n^{-1}\mathbf{P}_{\mathcal{B}_J}(\bphi(\bX_j)-\bphi(\bX))),
\end{eqnarray*}\fi
$K_m(\cdot)$ is a multivariate kernel function defined on $\mathds{R}^m$, e.g., $K_m(\bX)=\prod\limits_{i=1}^m K(\bX^{(i)})$ with $\bX^{(i)}$ being the $i$-th  element of $\bX$, and $h_n$ is a positive number depending on $n$, called the bandwidth or window width. By the estimation of~$\Lambda$, we have ~\begin{equation}\label{www}\widehat{w}(t',t,\mathbf{P}_{\mathcal{B}_J}\bphi(\bX))=\mathrm{exp}\{\widehat{\Lambda}(t',t|\mathbf{P}_{\mathcal{B}_J}\bphi(\bX))\}.\end{equation}
Then we can solve the generalized eigenvalue-eigenvector problem \eqref{eq10} to get $\left\{\widehat{\balpha}_j,\ j=1,2,\cdots,q\right\}$ and the estimator of~$u_j(\cdot)$ is
$$\widehat{u}_j(\bX)=\langle\widehat{\bm{\beta}}_j,\bphi(\bX)\rangle ={{R}}(\bX_{\le n},\bX)^T\widehat{\bm{\alpha}}_{j},\quad j=1,\cdots,q,$$
where~$\widehat{\bm{\beta}}_{j}=\Phi\widehat{\bm{\alpha}}_{j}$.

{
\section{Implementation of the proposed method}\label{tuningp}
In the following, we present an algorithm for implementing the proposed method.

\begin{enumerate}
  \item[Step 1] Chose a reproducing kernel $R$ and compute the Gram matrix ${{R}}$.

\item[Step 2] First apply double slicing on $(\widetilde{T},\Delta)$ to compute ${{Q}}_J$.
Then solve $${{R}} {{Q}}_J{{R}}\balpha=\lambda_J ({{R}}^2+n^2s {I}_n)\balpha$$ to obtain
the first $m$ significant eigenvectors $\widehat{\balpha}_{J1}^0$, $\ldots$, and $\widehat{\balpha}_{Jm}^0$,
and unitize them to $\widehat{A
}_J=(\widehat{\balpha}_{J1},\cdots,\widehat{\balpha}_{Jm})$ where
$$\widehat{\balpha}_{Jj}=\widehat{\balpha}_{Jj}^0\Big/\sqrt{\widehat{\balpha}_{Jj}^{0T}{{R}}\widehat{\balpha}_{Jj}^{0}}$$
satisfying that  $\langle\widehat{\bbeta}_{Jj},\widehat{\bbeta}_{Jj}\rangle=\langle\Phi\widehat{\balpha}_{Jj},\Phi\widehat{\balpha}_{Jj}\rangle=
\widehat{\balpha}_{Jj}^T{{R}}\widehat{\balpha}_{Jj}=1$.

\item[Step 3] First apply the $m$-dimensional kernel smoothing to obtain the weight function $\widehat{w}$, and then find $\widehat{W},\ \widehat{P}$ and ${{Q}}$. Afterwards, solve
$${{R}} {{Q}}{{R}}\balpha=\lambda ({{R}}^2+n^2\tau {I}_n)\balpha$$ to obtain
the first $q$ significant eigenvectors $\widehat{\balpha}_1$, $\ldots$, and $\widehat{\balpha}_q$ and then compute $$\widehat{u}_j(\bX)={{\bm R}}(\bX_{\le n},\bX)^T\widehat{\bm{\alpha}}_{j},\quad j=1,\cdots,q.$$
\end{enumerate}

In our simulation study in section \ref{ssimu}, we set the $L=L_0=L_1=10$ that is a moderate value for the number of slicing \cite{Li1991} and we show that the convergence rate is not affected by this number in section \ref{theory}. It is noted that $m$ and $q$, the numbers of significant eigenvectors, are decided in the same way as in principal component analysis (PCA). Here, we take the eigenvectors corresponding to the first $m$ (or $q$) eigenvalues that capture $90\%$ of the total sum of the eigenvalues.

It is important to choose a proper regularization parameter $\tau$ or $s$ as it controls the tradeoff between the bias and
the variance of the estimator. Similar to \cite{RSIR1991}, we consider the mean squared error (MSE) of the SDR directions
\begin{equation}\label{tauselect}
L(\tau)=V(\tau)+B(\tau)=\sum_{i=1}^n\sum_{j=1}^q\text{var}\Big(\widehat{u}_j(\bm{X}_i)(\tau)\Big)+\sum_{i=1}^{n}\sum_{j=1}^q
\Big(E\big(\widehat{u}_j(\bm{X}_i)(\tau)\big)-u_j(\bm{X}_i)\Big)^2,
\end{equation}
where $V(\tau)$ is the sum of the variances of $\widehat{u}_j(\bm{X_i})(\tau)$, $i=1,\ldots, n$, $j=1,\ldots,q$,
$B(\tau)$ is the corresponding sum of the squared biases. Here, we use the notation $\widehat{u}_j(\bm{X}_i)(\tau)$ instead of $\widehat{u}_j(\bm{X}_i)$ to emphasize the dependence of $\widehat{u}_j(\bm{X}_i)$ on $\tau$.
Let $\tau^*$ minimize $L(\tau)$. However, $L(\tau)$ is not directly computed since both $V(\tau)$ and $B(\tau)$ cannot be computed directly. To solve this problem, we can turn to a bootstrap procedure.
This strategy is motivated by \cite{RSIR1991}.
Specifically, we denote $\widehat{u}^{(b)}_1(\bm{X})(\tau),\ldots,\widehat{u}^{(b)}_q(\bm{X})(\tau)$ obtained from the $b$-th bootstrap sample, $b=1,\ldots,B$.  The estimates of $E(\widehat{u}_j(\bm{X}_i)(\tau))$ and $\text{var}(\widehat{u}_j(\bm{X}_i)(\tau))$ are thus given by
\begin{eqnarray}\label{bootstrap}
  \widehat{E}(\widehat{u}_j(\bm{X}_i)(\tau))&=&\dfrac{1}{B}\sum_{b=1}^{B}\widehat{u}^{(b)}_j(\bm{X}_i)(\tau)\\
  \widehat{\text{var}}(\widehat{u}_j(\bm{X}_i)(\tau))&=&\dfrac{1}{B-1}\sum_{b=1}^{B}\Big(\widehat{u}^{(b)}_j(\bm{X}_i)
  (\tau)-\widehat{E}(\widehat{u}_j(\bm{X}_i)(\tau))\Big)^2.
\end{eqnarray}
In practice, the true direction $u_j(\bm{X}_i)$ is also not unknown but it can be approximated by the bootstrap mean of
$E(\widehat{u}_j(\bm{X}_i)(\tau_0))$, where $\tau_0$ is a suitably chosen small positive value. Similar to \cite{functional2017},
we recommend to choose $\tau_0$ such that $n^2\tau_0=0.05{\lambda}_1({{R}}^2)$, where ${\lambda}_1(R^2)$ denotes the largest eigenvalue of ${{R}}^2$. We minimize $L(\tau)$ over a grid of $\tau$ in $[\tau_0/20,20\tau_0]$. The grid consists of 20 points, equally spaced in log scale. Although the proposed procedure for determining $\tau$ is based on an approximation, our experience from the extensive simulation study suggests that it works well.}

\section{Theoretical results}\label{theory}

In this section, we introduce four lemmas and present two theorems.

\subsection{Theoretical justifications of the proposed methodology}

In this subsection, we present three lemmas that provide theoretical justifications of the proposed methodology.

\begin{lemma}\label{lemma3.1}
{\it If for model (\ref{nonModle}), the linear design condition on the nonlinear dimension reduction holds, that is for $\forall \bm{\beta}\in \mathcal{H}$ and $E(\langle \bm{\beta},\bphi(\bX)\rangle|\langle\bm{\beta}_1,\bphi(\bX)\rangle,\cdots,\langle\bm{\beta}_q,\bphi(\bX)\rangle)$, there exits $c_0,c_1,\cdots,c_q$
such that~
\begin{equation}E(\langle \bm{\beta},\bphi(\bX)\rangle|\langle\bm{\beta}_1,\bphi(\bX)\rangle,\cdots,\langle\bm{\beta}_q,\bphi(\bX)\rangle)=c_0+c_1\langle\bm{\beta}_1,\bphi(\bX)\rangle+\cdots+c_q\langle\bm{\beta}_q,\bphi(\bX)\rangle,\label{lem1}\end{equation}
then $E(\bphi(\bX)|T)-E(\bphi(\bX))\in {\rm span}\{\Sigma\bm{\beta}_j,~j=1,\cdots,q\}.$}
\end{lemma}

The proof of Lemma \ref{lemma3.1} is given in \cite{Wu2008}.

The linear design condition above is a sufficient condition for kernel sliced inverse regression while the
performance of SIR or KSIR is not very sensitive to this condition as discussed by Li \cite{Li1991}. We can see
that the condition is fulfilled if the distribution of $\bphi(\bm{X})$ is elliptically symmetric such as normal distributions. In fact the low-dimensional projection of high-dimensional data often looks like normally distributed \cite{Diaconis1984}.
As a result , when the condition is mildly violated, the inverse regression procedure usually still works.

\begin{lemma}\label{equasolve}
{\it~ Solutions of the sample version \eqref{genalesti} of the infinite dimensional generalized eigenvalue-eigenvector problem \eqref{eq7} are as follows:
\begin{eqnarray}\label{lemma2.01}
\widehat{\bm{\beta}}_{j}=\sum_{i=1}^n\widehat{\bm{\alpha}}_{j}^{(i)}\bphi(\bX_i)=\Phi\widehat{\bm{\alpha}}_{j},\quad j=1,\cdots,q,
\end{eqnarray}
where~$\widehat{\bm{\alpha}}_{1},\cdots,\widehat{\bm{\alpha}}_{q}$ are the solutions of the following dimensional generalized eigenvalue-eigenvector problem:
\begin{eqnarray}\label{lemma2.01b}
{{R}}{{Q}}{{R}}{\bm{\alpha}}={\lambda}{{R}}^2{\bm{\alpha}},
\end{eqnarray}
which is equivalent to \eqref{genalesti}. Similar results hold true for using the following regularization technique:
\begin{eqnarray}\label{lemma2.02}
\widehat{\Sigma}\widehat{\Gamma}\bm{\beta} = \lambda(\widehat{\Sigma}^2+\tau I_n)\bm{\beta},
\end{eqnarray}
i.e., the solutions of \eqref{lemma2.02} have the same type of expressions as~\eqref{lemma2.01} with ~${\widehat\balpha}_j$s being the solutions of ~${{R}}{{Q}}{{R}}{\bm{\alpha}}={\lambda}({{R}}^2+n^2\tau{I}_n){\bm{\alpha}}$.}
\end{lemma}

The proof of Lemma \ref{equasolve} is given in the appendix A.

Results of Lemma~\ref{equasolve} can be extended to the case corresponding to joint~SDR. When~${{Q}}$ is replaced by  ${{Q}}_J$, similar results can be obtained.
Proof of Lemma~\ref{equasolve} is similar to the property 7 in \cite{Wuqiang2013}. In stead of considering complete data as in \cite{Wuqiang2013},
we deal with censored survival data and consider the information of censoring in ${{Q}}$ in this paper. The result \eqref{lemma2.01} is of great importance since it converts an infinite dimensional problem to a finite one and hence the proposed methodology is applicable.

\begin{lemma}\label{conditionsurv}
{\it If the condition 1 holds, then we have
  \begin{eqnarray*}
   w(t',t,\bphi(\bX))&=&\mathrm{exp}\{-\Lambda(t',t|\bphi(\bX))\}\\
    E\big[\bphi(\bX)I(T\geq t)\big]&=& E\big[\bphi(\bX)I(\widetilde{T}\geq t)\big]+E\big[\bphi(\bX)I(\widetilde{T}< t, \Delta=0)w(\widetilde{T},t,\bphi(\bX))\big], \\
    E\big[I(T\geq t)\big]&= &E\big[I(\widetilde{T}\geq t)\big]+E\big[I(\widetilde{T}< t, \Delta=0)w(\widetilde{T},t,\bphi(\bX))\big],
  \end{eqnarray*}
where
\begin{equation}\label{wu}
  \Lambda(t',t|\bphi(\bX))=E\left\{\frac{I(t'\leq\widetilde{T}<t,\Delta=1)}{S_{\widetilde{T}}(\widetilde{T}|\bphi(\bX))}\Big|\bphi(\bX)\right\}.
\end{equation}
}
\end{lemma}

The proof of Lemma \ref{conditionsurv} is given in the appendix B.

We know that $T$ is unobservable true lifetime, and hence we cannot estimate $E[\bphi(\bm{X})I(T\geq t)]$ and $E[I(T\geq t)]$ directly. Lemma \ref{conditionsurv} provides a way to estimate them through the conditional survival function $w(t',t,\bphi(\bX))$. Therefore, the remaining problem is to estimate $w(t',t,\bphi(\bX))$, or equivalently, to estimate the conditional cumulative hazard function $\Lambda(t',t|\bphi(\bm{X}))$. By \eqref{wu}, $\Lambda(t',t|\bphi(\bm{X}))$ is a conditional expectation conditioned on $\bphi(\bm{X})$.  In terms of  \eqref{omegaPb},  $\bphi(\bm{X})$ can be replaced by $\mathbf{P}_{\mathcal{B}_J}\bphi(\bX)$ in \eqref{wu}, i.e., $\Lambda(t',t|\bphi(\bm{X}))=\Lambda(t',t|\mathbf{P}_{\mathcal{B}_J}\bphi(\bX))$, which can be estimated by a double slicing procedure.

\subsection{Asymptotic properties}

In this subsection, we study the asymptotic properties of the sufficient dimension directions estimated by RDSIR. We first need to show that the following two conditions hold true:

\qquad$(L1)~\quad\left|{\widehat{p}_\ell}-p_\ell\right|= O_p(n^{-\kappa}),$

\qquad$(L2)~\quad\left\|\dfrac{1}{n}\Phi\widehat{\bm{w}}_\ell-E\Big[\bphi(\bX)I(T\in D_\ell)\Big]\right\|=O_p(n^{-\kappa}),$

\noindent where $\kappa\in (0,1/4)$. We note that $$\widehat{p}_l=\widehat{P}(T\geq t_l)-\widehat{P}(T\geq t_{l+1})=\dfrac{1}{n}\mathbf{1}_n^T\widehat{\bm{w}}_{\ell}.$$
To prove $(L1)$, we need only to show
\begin{eqnarray}
\Big|\widehat{P}(T\geq t)-P(T\geq t)\Big|=O_p(n^{-\kappa})
\end{eqnarray}
where
\begin{eqnarray}\label{phat}
\widehat{P}(T\geq t)=\dfrac{1}{n}\sum\limits_{i=1}^{n}I(\widetilde{T}_i\geq t)+\dfrac{1}{n}\sum\limits_{i:\widetilde{T}_i<t,\triangle_i=0}\widehat{w}\big(\widetilde{T}_i,t,\bphi(\bm{X}_i)\big).
\end{eqnarray}

As \eqref{omegaPb} implies that $w(t',t,\bphi(\bm{X}_i))=w(t',t,\bm{P}_{\mathcal{B}_J}\bphi(\bm{X}_i))$, we can substitute
$\bphi(\bm{X}_i)$ in \eqref{phat} with $\bm{P}_{\mathcal{B}_J}\bphi(\bm{X}_i)$. We denote a basis of $\mathcal{B}_J$ by
$B_J=\{\bbeta_{J1},\cdots,\bbeta_{Jm}\}$, which is found via the generalized eigenvalue-eigenvector problem \eqref{eq15}, and its corresponding estimator by $\widehat{B}_J=\{\widehat{\bbeta}_{J1},\cdots,\widehat{\bbeta}_{Jm}\}$ that is obtained via the regularization problem \eqref{rr}.

To simplify notations, we denote $\bm{Z}_i=B_J^T\bm{\phi}(\bm{X}_i)$ and $\widehat{\bm{Z}}_i=\widehat{B}_J^T\bm{\phi}(\bm{X}_i)$, where $\bm{Z}_i$ is an $m$-dimensional vector whose $k$-th component is $\bm{Z}_i^{(k)}=\langle \bbeta_{Jk},\bm{\phi}(\bm{X}_i) \rangle$, and $\widehat{\bm{Z}}_i$ is the corresponding estimator of $\bm{Z}_i$.
To further simplify the notations, we put
\begin{align}\label{notations}
      e_{ij} & =h_n^{-m}K_m\big(h_n^{-1}(\bm{Z}_i-\bm{Z}_j)\big)\ ,\quad
      u_{ij} =e_{ij}-E(e_{ij}|\bm{Z}_j).\notag \\
      v_{kj} &=I(\widetilde{T}_k>\widetilde{T}_j)e_{kj}-
E\left[I(\widetilde{T}_k>\widetilde{T}_j)e_{kj}\Big|\bm{Z}_j,\widetilde{T}_j\right]\notag\\
      f_j &=f(\bm{Z}_j)\ ,\quad S_j=S_{\widetilde{T}}(\widetilde{T}_j|\bm{Z}_j)\ ,\quad I_{ij}=I(\widetilde{T}_i<\widetilde{T}_j<t,\triangle_j=1)\notag\\
      \Lambda_i &=\Lambda(\widetilde{T}_i,t|\bm{Z}_i)\ ,\quad \widehat{\Lambda}_i=\widehat{\Lambda}(\widetilde{T}_i,t|\bm{Z}_i)\ ,\quad w_i=e^{-\Lambda_i}.
    \end{align}

We make the following assumptions:
\renewcommand{\theequation}{A\arabic{equation}}
\setcounter{equation}{0}
\begin{flalign}\label{A1}
  &\qquad (\text{A1})\ E(e_{ij}|\bm{Z}_j)-f(\bm{Z}_j)=O_p(n^{-1/2}).&
\end{flalign}
\vspace{-1cm}
\begin{flalign}\label{A2}
  &\qquad(\text{A2})\ \dfrac{1}{n}\sum\limits_{k}u_{ki}=O_p(n^{-1/4}).&
\end{flalign}
\vspace{-1cm}
\begin{flalign}\label{A3}
  &\qquad(\text{A3})\ E\left[I(\widetilde{T}_k>\widetilde{T}_j)e_{kj}|\bm{Z}_j,\widetilde{T}_j\right]
-S_{\widetilde{T}}(\widetilde{T}_j|\bm{Z}_j)f(\bm{Z}_j)=O_p(n^{-1/2}).&
\end{flalign}
\vspace{-1cm}
\begin{flalign}\label{A4}
  &\qquad(\text{A4})\ \dfrac{1}{n}\sum\limits_{k}v_{ki}=O_p(n^{-1/4}).&
\end{flalign}
\vspace{-1cm}
\begin{flalign}\label{A5}
  &\qquad(\text{A5})\ E\left(I_{ij}S_j^{-1}e_{ji}|\bm{Z}_i,\widetilde{T}_i\right)-\Lambda_if_i=O_p(n^{-1/2}).&
\end{flalign}
\vspace{-1cm}
\begin{flalign}\label{A6}
  &\qquad(\text{A6})\ \int x^iK(x)dx=0,\ i=1,\cdots,d-1\ \text{and}\ \int x^dK(x)dx\ne 0\ \text{where}\ d>(m+1)/b&\nonumber\\ &\qquad\qquad\quad\text{with}\ 0<b\leq 1/4.&
\end{flalign}
\vspace{-1cm}
\begin{flalign}\label{A7}
  &\qquad(\text{A7})\ \text{The tuning parameter}\ s=s(n)\ \text{in} \ \eqref{rr}\ \text{satisfies}\ \lim\limits_{n\rightarrow \infty}s=0\ \text{and}\ \lim\limits_{n\rightarrow \infty}s\sqrt{n}=\infty.&
\end{flalign}
\vspace{-1cm}
\begin{flalign}\label{A8}
  &\qquad(\text{A8})\ \{\bbeta_{Jj}\}_{j=1}^m\ \text{depends only on finite eigenfunctions of}\ \Sigma.&
\end{flalign}
\vspace{-1cm}
%\begin{flalign}\label{A9}
%  &\qquad(\text{A9})\ \mathscr{R}(\bm{s},\bm{t})\ \text{is a continuous function and the value of}\ \bm{X}\ \text{is taken from a compact set}.&
%\end{flalign}
\renewcommand{\theequation}{\arabic{equation}}
\setcounter{equation}{39}

It can be observed that the assumptions \eqref{A1}-\eqref{A5} are the regularity conditions that are the parallel extensions  of the assumptions made in Lemma 3.1 of \cite{Li1999} and are satisfied with
bandwidth $h_n\propto n^{-1/2d}$ provided $m\leq d$.
The assumption \eqref{A1} implies that the bias term of $\widehat{f}(\bm{Z}_i)$ is of root $n$ rate which requires a bandwidth smaller than the usual optimal one. The assumption \eqref{A2} is a mild and flexible one which only requires that the rate of the term $\sum\limits_{k}u_{ki}/n$ contributing to the variance of $\widehat{f}(\bm{Z}_i)$ is $O_p(n^{-1/4})$.
The assumptions \eqref{A3} and \eqref{A5} are made so that the biases of both kernel estimates of $S_jf_j$ and $\Lambda_if_i$ have the root $n$ rate. With suitable smoothness conditions on conditional survival function $S_{\widetilde{T}}(t|\bm{Z})$ and cumulative hazard function $\Lambda(t',t|\bm{Z})$, the bandwidth to achieve the assumption \eqref{A1} may also imply that both assumptions \eqref{A3} and \eqref{A5} hold true. The assumption \eqref{A4} is made for the same reason as the assumption \eqref{A2}.  The assumption \eqref{A6} is imposed on the kernel function, which is similar to those made in \cite{Li1999}, which also ensures the rationality of the results. It is noted that a relatively large $d$ can be chosen such that $\kappa=b-(m+1)/d>0$. We can refer to \cite{Jones1993} for the construction of this kind of higher order kernel.
The assumptions \eqref{A7}-\eqref{A8} that are the same assumptions in Theorem 9 of \cite{Wuqiang2013}, which have been used to study the convergence rates of the joint SDR directions for both survival and censoring times there. Under the assumptions \eqref{A7}-\eqref{A8}, we can learn from the results of
Wu \textit{et al.} \cite{Wuqiang2013} that $\widehat{\bm{Z}}_i^{(k)}-\bm{Z}_i^{(k)}=O_p(n^{-b})$ where $0<b\leq 1/4$.
%The assumption \eqref{A9} is a standard and common condition for the reproducing kernel and we impose
%the value of $\bm{X}$ is taken from a compact set is also reasonable in practical application. Besides, under assumption \eqref{A9}, $\int \mathscr{R}(\bm{X},\bm{X})d\bm{X}=\sum_ia_i<\infty$, which is important in the proof.

We can now formally introduce the following important lemma.

\begin{lemma}\label{lemma4}
Under the assumptions \eqref{A1}-\eqref{A8}, there exits a positive number $\kappa\in (0,1/4)$ such that

\quad$(L1)~\left|{\widehat{p}_\ell}-p_\ell\right|= O_p(n^{-\kappa})$ %where $\widehat{p}_\ell$ is given in

\quad$(L2)~\left\|\dfrac{1}{n}\Phi\widehat{\bm{w}}_\ell-E\Big[\bphi(\bX)I(T\in D_\ell)\Big]\right\|=O_p(n^{-\kappa}).$
\end{lemma}

The proof of Lemma \ref{lemma4} is given in the appendix C.

Our asymptotic results are based on the perturbation theory for linear operators. We first introduce some properties of Hilbert-Schmidt operators. Let $\mathcal{H}$ be a Hilbert space and $\{\psi_i, i\in I\}$ a standard orthogonal basis of $\mathcal{H}$. We say that linear operator $\mathcal{L}$ defined on~$\mathcal{H}$ is a Hilbert-Schmidt operator if
\begin{equation*}
  \|\mathcal{L}\|_{HS}^2=\sum_{i=1}^\infty \|\mathcal{L}\psi_i\|_{\mathcal{H}}^2<\infty.
\end{equation*}
Hilbert-Schmidt class spans as a new~Hilbert space with~$\|\cdot\|_{HS}$ norm.
If~$\mathcal{S}$ is a bounded operator in $\mathcal{H}$, then $\mathcal{S}\mathcal{L}$ and~$\mathcal{L}\mathcal{S}$ belongs to~Hilbert-Schmidt class and
\begin{equation*}
  \|\mathcal{S}\mathcal{L}\|_{HS}\leq \|\mathcal{S}\|\|\mathcal{L}\|_{HS}, ~~ \|\mathcal{L}\mathcal{S}\|_{HS}\leq \|\mathcal{L}\|\|\mathcal{S}\|_{HS},
\end{equation*}
where~$\|\cdot\|$ is the norm of operator
\begin{equation*}
 \|\mathcal{L}\|=\sup_{f\in\mathcal{H}}\frac{\|\mathcal{L}f\|}{\|f\|}.
\end{equation*}
The next two theorems are about the
convergence rates of the SDR directions corresponding to the true survival time.

\begin{theorem}\label{limit01}
{Under the assumptions of \emph{Lemma \ref{lemma4}}, for any $N>0$, we have
\begin{equation}\label{d3.5}
  \left\|\left({\widehat{\Sigma}}^2+\tau I\right)^{-1}\widehat{\Sigma}\widehat{\Gamma}-\Sigma^{-1}\Gamma\right\|_{HS}=O_p\left(\frac{1}{\tau n^{\kappa}}\right)+\sum_{j=1}^{d_\Gamma}\left(\frac{\tau}{\nu_N^2}\Big\|\Psi_N(\widetilde{\varphi}_j)\Big\|+\Big\|\Psi_{N}^{\perp}(\widetilde{\varphi}_j)\Big\|\right)
\end{equation}
where~$\widetilde{\varphi}_j=\Sigma^{-1}\varphi_j$, $\Psi_N$ and $\Psi_{N}^{\perp}$ are respectively the projection operator and its complement,
$$\Psi_{N}=\sum_{j=1}^N\psi_j\otimes\psi_j,\quad \Psi_{N}^{\perp}=\sum_{j=N+1}^\infty\psi_j\otimes\psi_j.$$
If the smoothing parameter $\tau=\tau(n)$ satisfies that $\lim\limits_{n\rightarrow\infty}\tau=0$ and $\lim\limits_{n\rightarrow\infty}\tau n^{\kappa}=\infty$, then
\begin{equation}\label{d3.6}
  \left\|({\widehat{\Sigma}}^2+\tau I)^{-1}\widehat{\Sigma}\widehat{\Gamma}-\Sigma^{-1}\Gamma\right\|_{HS}=o_p(1).
\end{equation}
}
\end{theorem}

The proof of Theorem \ref{limit01} is given in the appendix D.

\begin{theorem}\label{limit}

{If $\lim\limits_{n\rightarrow\infty}\tau=0, \lim\limits_{n\rightarrow\infty}\tau n^{\kappa}=\infty$, $d_\Gamma$ is the rank of~~${\Gamma}$, ~$\{\widehat{\bm{\beta}}_j\}$ is the estimation from~(\ref{genalesti}), then
\begin{equation}\label{d3.6bbb}
  \big|\widehat{u}_j(\bX)-u_j(\bX)\big|=\left|\langle \widehat{\bm{\beta}}_j,\bphi(\bX)\rangle-\langle {\bm{\beta}}_j,\bphi(\bX)\rangle\right|=o_p(1),~~j=1,\cdots,d_\Gamma,
\end{equation}
further more, under assumption \eqref{A8}, $\{\bm{\beta}_j\}_{j=1}^{d_\Gamma}$ depend only on a finite number of eigenvectors of the covariance operator of $\Sigma$, the rate of convergence is~$O(n^{-\kappa/2})$. That is
}
\begin{equation}\label{d3.6bbbb}
  \big|\widehat{u}_j(\bX)-u_j(\bX)\big|=\left|\langle \widehat{\bm{\beta}}_j,\bphi(\bX)\rangle-\langle {\bm{\beta}}_j,\bphi(\bX)\rangle\right|=O_p(n^{-\kappa/2}),~~j=1,\cdots,d_\Gamma.
\end{equation}

\end{theorem}

Theorem \ref{limit} is a corollary of Theorem \ref{limit01}. By a well-known result in the perturbation theory \cite{Kato2013}, the eigenspaces of $(\widehat{\Sigma}^2+\tau I)^{-1}\widehat{\Sigma}\widehat{\Gamma}$ converge to those of $\Sigma^{-1}\Gamma$ at the same rate. Therefore, the proof of Theorem \ref{limit} is omitted.

\section{Simulations and real data analysis}\label{sec:Simu}
\subsection{Simulations}\label{ssimu}
In this section, we carry out simulation studies of RDSIR, and we also compare RDSIR with double sliced
inverse regression (DSIR) proposed by (\cite{Li1999}) which is suitable for linear sufficient dimension reduction.

Commonly used kernels are the Gaussian radial basis kernel ~$\mathscr{R}(s,t)=\mathrm{exp}(-\mathrm{scale}\cdot\| s-t\|^2)$ and the polynomial kernel~$\mathscr{R}(s,t)=(\mathrm{scale}\cdot\langle s,t\rangle+\mathrm{offset})^{\mathrm{degree}}$. \cite{Duan2003,Keerthi2003} proposed some criterion about how to choose reproducing kernel and corresponding parameter.  When there is no prior information on a dataset, the Gaussian radial basis kernel is usually chosen.
In the following examples, we set training sample size~$n_{tr}=100$ and test sample size~$n_{te}=200$. The regularization parameters $\tau$ and $\tau_J$ can be chosen by the criterion \eqref{tauselect} proposed in the end of the section \ref{ecs}.

To compare the estimated and true dimension reduction directions, which may be vectors of different dimensions, we adopt the \textit{Robust Maximum Association Estimators} (\textbf{RMAE}) (\cite{Alfons2017}) as follows. Suppose that $\bm{X}$ is a $p$-dimensional random vector and $\bm{Y}$ is a $q$-dimensional random vector, with $p\geq q$. A measure of multivariate association between $\bm{X}$ and $\bm{Y}$ can be defined by looking for linear combinations $\bm{\alpha}^T\bm{X}$ and $\bm{\beta}^T\bm{Y}$ of
the original variables that has the maximal association. That is, we seek a measure
\begin{equation}\label{}
 \rho_r(\bX,\bY)=\max_{\balpha,\bbeta}\ r(\balpha^T\bm{X},\bbeta^T\bm{Y}),
\end{equation}
where $r$ is a measure of association between two univariate variables. Taking the classical Pearson correlation
for $r$ results in the first canonical correlation coefficient (\cite{Johnson02}). The bivariate association measure $r$ considered in this article is Spearman's rank correlation which is defined as
\begin{equation}\label{}
 r(U,V)=r_P(\text{rank}(U),\text{rank}(V)),
\end{equation}
where $r_P(X,Y)=\text{cov}(X,Y)/\sqrt{\text{var}(X)\text{var}(Y)}$ is the Pearson correlation and\ $\text{rank}(u)=F_W(u)$, with $F_W$ the cumulative distribution function of the random variable $W$. For practical application, an empirical version of $\rho_r$ denoted by $\hat{\rho}_r$ is needed. \cite{Alfons2017} developed the \emph{alternate grid} algorithm for the computation of such maximum association estimates and studied their theoretical properties for various association measures such as Pearson, Spearman and Kendall's $\tau$ correlation. It turns out that the Spearman and Kendall's $\tau$ correlation
yield a maximum association estimate with good robustness properties and good efficiency, which has been implemented in R package {\tt ccaPP} (\cite{Alfons2016}).
In the following, we perform 100 Monte Carlo simulations for each model and report the averages and standard deviations of \textbf{RMAE} under different censoring proportions.

\begin{table}[t]
\centering\caption{Mean and standard deviation of the \textbf{RMAE} under different censoring proportion}\label{example12}
\begin{tabular}{llcllll}
\hline
                         & \multicolumn{1}{c}{}                     &                           & \multicolumn{1}{c}{0\%}          & \multicolumn{1}{c}{20\%}         & \multicolumn{1}{c}{40\%}         & \multicolumn{1}{c}{60\%}         \\ \cline{3-7}
\multirow{4}{*}{Model 1} & \multicolumn{1}{c}{\multirow{2}{*}{$q=1$}} & RDSIR                     & \multicolumn{1}{c}{0.933(0.027)} & \multicolumn{1}{c}{0.883(0.044)} & \multicolumn{1}{c}{0.801(0.063)} & \multicolumn{1}{c}{0.728(0.091)} \\
                         & \multicolumn{1}{c}{}                     & DSIR                      & \multicolumn{1}{c}{0.936(0.024)} & \multicolumn{1}{c}{0.895(0.036)} & \multicolumn{1}{c}{0.828(0.055)} & \multicolumn{1}{c}{0.766(0.076)} \\ \cline{3-7}
                         & \multirow{2}{*}{$q=2$}                     & RDSIR                     & 0.939(0.024)                     & 0.896(0.036)                     & 0.823(0.058)                     & 0.762(0.070)                     \\
                         &                                          & DISR                      & 0.940(0.022)                     & 0.903(0.032)                     & 0.843(0.053)                     & 0.789(0.064)                     \\ \cline{2-7}
\multirow{4}{*}{Model 2} & \multirow{2}{*}{$q=1$}                     & \multicolumn{1}{l}{RDSIR} & 0.895(0.036)                     & 0.886(0.035)                     & 0.838(0.050)                     & 0.792(0.059)                     \\
                         &                                          & \multicolumn{1}{l}{DSIR}  & 0.903(0.028)                     & 0.899(0.030)                     & 0.867(0.039)                     & 0.834(0.045)                     \\ \cline{3-7}
                         & \multirow{2}{*}{$q=2$}                     & \multicolumn{1}{l}{RDSIR} & 0.903(0.033)                     & 0.897(0.031)                     & 0.854(0.046)                     & 0.814(0.057)                     \\
                         &                                          & \multicolumn{1}{l}{DSIR}  & 0.910(0.027)                     & 0.907(0.029)                     & 0.878(0.034)                     & 0.849(0.042)                     \\ \hline
\end{tabular}
\end{table}

\noindent\textbf{Model 1 (Linear)}. We generate survival time and censoring time from
\begin{equation}\label{example2sim}
\begin{split}
& T= \Phi\big(\epsilon*\bbeta_1^T\bX\big) \\
& C= \Phi(\bbeta_2^T\bX)+U
\end{split}
\end{equation}
where $\bX\sim N(0,I_{50})$, $\bbeta_1=(1,0,-1,0,\cdots,1,0,-1,0,0,0,\cdots,0)^T/5$, the last 10 elements are 0. $\bm{\beta}_2=(0,0,\cdots,0,1,-1,\cdots,1,-1)^T/5$, the first 10 elements are 0. $\Phi$ is the cumulative distribution function of the standard normal distribution and $\epsilon$ follows an uniform distribution $U(0,1)$. The symbol $U$ denotes a random variable uniformly distributed on $(0,c)$, where $c$ controls the censoring proportion. We set the reproducing kernel $\mathscr{R}(s,t)=\left\langle s,t\right\rangle$ and consider the first two directions. The results are given in Table \ref{example12}.

\begin{table}[t]
\centering\caption{The \textbf{RMAE} mean and standard deviation under different censoring proportion}\label{example3}
\begin{tabular}{ccccccc}
\hline
                                                                              &                      &       & 0\%          & 20\%         & 40\%         & 60\%         \\ \cline{3-7}
\multirow{4}{*}{\begin{tabular}[c]{@{}c@{}}Model 3\\ $(n=100,p=50)$\end{tabular}} & \multirow{2}{*}{$q=1$} & RDSIR & 0.966(0.033) & 0.928(0.107) & 0.898(0.186) & 0.866(0.213) \\
                                                                              &                      & DSIR  & 0.065(0.051) & 0.055(0.038) & 0.068(0.048) & 0.062(0.043) \\ \cline{3-7}
                                                                              & \multirow{2}{*}{$q=2$} & RDSIR & 0.972(0.015) & 0.948(0.075) & 0.934(0.120) & 0.929(0.099) \\
                                                                              &                      & DISR  & 0.103(0.050) & 0.086(0.041) & 0.106(0.044) & 0.100(0.043) \\ \hline
\multirow{4}{*}{\begin{tabular}[c]{@{}c@{}}Model 3\\ $(n=100,p=60)$\end{tabular}} & \multirow{2}{*}{$q=1$} & RDSIR & 0.873(0.024) & 0.837(0.115) & 0.814(0.128) & 0.787(0.168) \\
                                                                              &                      & DSIR  & 0.055(0.042) & 0.051(0.046) & 0.058(0.050) & 0.051(0.045) \\ \cline{3-7}
                                                                              & \multirow{2}{*}{$q=2$} & RDSIR & 0.877(0.021) & 0.859(0.057) & 0.848(0.069) & 0.840(0.079) \\
                                                                              &                      & DSIR  & 0.097(0.044) & 0.089(0.048) & 0.094(0.052) & 0.093(0.048) \\ \hline
\multirow{4}{*}{\begin{tabular}[c]{@{}c@{}}Model 3\\ $(n=100,p=70)$\end{tabular}} & \multirow{2}{*}{$q=1$} & RDSIR & 0.783(0.092) & 0.757(0.122) & 0.724(0.147) & 0.695(0.169) \\
                                                                              &                      & DSIR  & 0.059(0.043) & 0.049(0.040) & 0.061(0.045) & 0.062(0.046) \\
                                                                              & \multirow{2}{*}{$q=2$} & RDSIR & 0.799(0.042) & 0.790(0.050) & 0.768(0.083) & 0.757(0.092) \\
                                                                              &                      & DSIR  & 0.103(0.042) & 0.090(0.040) & 0.099(0.043) & 0.093(0.044) \\ \hline
\end{tabular}
\end{table}
\noindent\textbf{Model 2 (Linear)}.
In this example, we consider two linear SDR directions $\bbeta_1^T\bX$ and $\bbeta_2^T\bX$. We simulate the survival time by
\begin{equation}\label{PHsim}
\begin{split}
 T =-\dfrac{\log(U)}{\exp(\bm{\beta}_1^T\bX)+\exp(\bm{\beta}_2^T\bX)}
\end{split}
\end{equation}
where $\bX\sim N(0,I_{50})$ and $U\sim U(0,1)$. The coefficient vector $\bm{\beta}_1=(1,-1,1,-1,\cdots,0,0)^T/5$ where the last 10 elements are 0 and $\bm{\beta}_2=(0,0,\cdots,1,-1,1,-1)^T/5$ where the first 10 elements are 0. The censoring time is generated from $U(0,c)$ and $c$ controls the censoring rate. It is obvious that this is a linear dimension reduction problem and we want to recover the space spanned by $\{\bm{\beta}_1,\bbeta_2\}$. We set the reproducing kernel $\mathscr{R}(s,t)=\left\langle s,t\right\rangle$
and the results are summarized in Table \ref{example12}.

It is noted that for both Model 1 and Model 2, the true SDR directions are in the 50-dimensional Euclidean space. By Table \ref{example12}, we can see that DSIR performs well on this type of linear dimension reduction problem. We can also learn from Table \ref{example12} that RDSIR can achieve similar performance compared to DSIR. In fact, the key idea of kernel trick is to obtain the dot
product in the high dimensional feature space by computing the value of the kernel function in the original input space. When using a linear kernel, RDSIR is equivalent to DSIR. There will be a small difference between RDSIR and DSIR due to the accumulation of calculation errors.

\begin{table}[t]
\centering\caption{Mean and standard deviation of the \textbf{RMAE} under different censoring proportion}\label{example4}
\begin{tabular}{ccccccc}
\hline
                                                                                &                      &       & 0\%          & 20\%         & 40\%         & 60\%         \\ \cline{3-7}
\multirow{4}{*}{\begin{tabular}[c]{@{}c@{}}Model 4\\ $(n=100,p=50)$\end{tabular}} & \multirow{2}{*}{$q=1$} & RDSIR & 0.895(0.125) & 0.884(0.176) & 0.873(0.178) & 0.851(0.187) \\
                                                                                &                      & DSIR  & 0.101(0.040) & 0.102(0.043) & 0.103(0.041) & 0.105(0.044) \\ \cline{3-7}
                                                                                & \multirow{2}{*}{$q=2$} & RDSIR & 0.934(0.057) & 0.932(0.092) & 0.930(0.067) & 0.901(0.134) \\
                                                                                &                      & DISR  & 0.136(0.039) & 0.147(0.043) & 0.141(0.043) & 0.140(0.046) \\ \hline
\multirow{4}{*}{\begin{tabular}[c]{@{}c@{}}model 4\\ $(n=100,p=60)$\end{tabular}} & \multirow{2}{*}{$q=1$} & RDSIR & 0.826(0.053) & 0.800(0.084) & 0.760(0.142) & 0.719(0.184) \\
                                                                                &                      & DSIR  & 0.092(0.044) & 0.106(0.060) & 0.082(0.039) & 0.108(0.057) \\ \cline{3-7}
                                                                                & \multirow{2}{*}{$q=2$} & RDSIR & 0.855(0.029) & 0.839(0.042) & 0.809(0.105) & 0.806(0.064) \\
                                                                                &                      & DSIR  & 0.141(0.050) & 0.138(0.045) & 0.130(0.044) & 0.135(0.067) \\ \hline
\multirow{4}{*}{\begin{tabular}[c]{@{}c@{}}Model 4\\ $(n=100,p=70)$\end{tabular}} & \multirow{2}{*}{$q=1$} & RDSIR & 0.758(0.050) & 0.721(0.078) & 0.708(0.118) & 0.631(0.224) \\
                                                                                &                      & DSIR  & 0.089(0.031) & 0.110(0.048) & 0.100(0.028) & 0.095(0.024) \\
                                                                                & \multirow{2}{*}{$q=2$} & RDSIR & 0.766(0.049) & 0.748(0.049) & 0.737(0.075) & 0.702(0.142) \\
                                                                                &                      & DSIR  & 0.144(0.035) & 0.150(0.040) & 0.140(0.039) & 0.127(0.040) \\ \hline
\end{tabular}
\end{table}
\noindent\textbf{Model 3 (Nonlinear)}. To show the performance of RDSIR when there exists a nonlinear dimension reduction subspace, we simulate the survival data by
\begin{equation}\label{NONAFTsim}
\begin{split}
&  T = U_1*\left[\bX^{(1)2}+\cdots+\bX^{(50)2}\right]\\
&  C = U_2*\left[\sin(X^{(1)})+\cdots+\sin(X^{(50)})\right]^2
\end{split}
\end{equation}
where $\bX\sim N(0,I_{p})$, $p=50,60,70$, $X^{(i)}$ is the $i$th element of $\bX$. The random error $U_1\sim U(0,1)$, $U_2\sim U(0,c)$ and $c$ controls the censoring proportion.
Let the reproducing kernel be $\mathscr{R}(s,t)=(\langle s,t\rangle+1)^2$, see Table \ref{example3} for the simulation results.

\noindent\textbf{Model 4 (Nonlinear)}. We consider another nonlinear dimension reduction example which exists two directions.
The survival data is simulated by
\begin{equation}\label{NONsim4}
  T = -\dfrac{\log(U)}{\sin^2(\bX^{(1)}/2)+\sin^2(\bX^{(2)}/2)+\cdots+\sin^2(\bX^{(50)}/2)}+\exp\left\{-\dfrac{1}{2}\bX^T\bX\right\},
\end{equation}
where the random variable $U\sim U(0,1)$
and censoring time $C\sim U(0,c)$ with $c$ controling the censoring rate, and $\bX\sim N(0,I_{p})$, $p=50,60,70$. We choose the Gaussian radial basis kernel $\mathscr{R}(s,t)=\mathrm{exp}(-\mathrm{scale}\cdot\| s-t\|^2)$. The results are displayed in  Table \ref{example4}. {
What we need to pay attention to is that both nonlinear SDR directions $u_1(\bX)=\sin^2(\bX^{(1)}/2)+\sin^2(\bX^{(2)}/2)+\cdots+\sin^2(\bX^{(50)}/2)$ and $u_2(\bX)=\bX^T\bX$ are artificially assumed in the model setting stage. We can also regard these two directions as one direction since they can be absorbed in the unknown function $g$. Thus, we may find that the first estimated direction contains most of the information and there is no significant improvement in the results after adding the second direction. On the other hand, it is also possible that the first two directions are not enough to characterize the relationship between $T$ and $\bX$. So the key problem is to choose a suitable dimension $q$ which has been discussed in Section \ref{tuningp}.}

By taking these examples into account, we can conclude that the larger the censoring rate and the number of regressors, the harder the problem and one might expect the performance of our method to deteriorate.
From Tables \ref{example12}-\ref{example4}, we can see that DSIR fails to extract nonlinear dimension reduction directions, but RDSIR proposed in this paper performs well in both linear and nonlinear situations.

\subsubsection{NCCTG Lung Cancer Data}

Consider the North Central Cancer Treatment Group (NCCTG) lung cancer data in \cite{lungdata}.
NCCTG lung cancer dataset recorded the survival days of patients with advanced lung cancer, together with assessments of the patients performance status measured either by the physician and by the patients themselves. There were total 228 patients that includes 63 patients whose data were right censored. Ten variates are included in our study, specifically, survival time in days (time: $T$), censoring status $1=$ censored and $2=$ dead (status: $\delta$), institution code (inst: $\bX_1$), age in years (age: $\bX_2$), male=1 and female=2 (sex: $\bX_3$),
ECOG performance score (0=good---5=dead) (ph.ecgo: $\bX_4$), Karnofsky performance score (bad=0---good=100) rated by physician (ph.karno: $\bX_5$), Karnofsky performance score as rated by patient (pat.karno: $\bX_6$), calories consumed at meals (meal.cal: $\bX_7$), weight loss in last six months (wt.loss: $\bX_8$).

\begin{figure}[h]
  \centering
  \includegraphics[width=8cm]{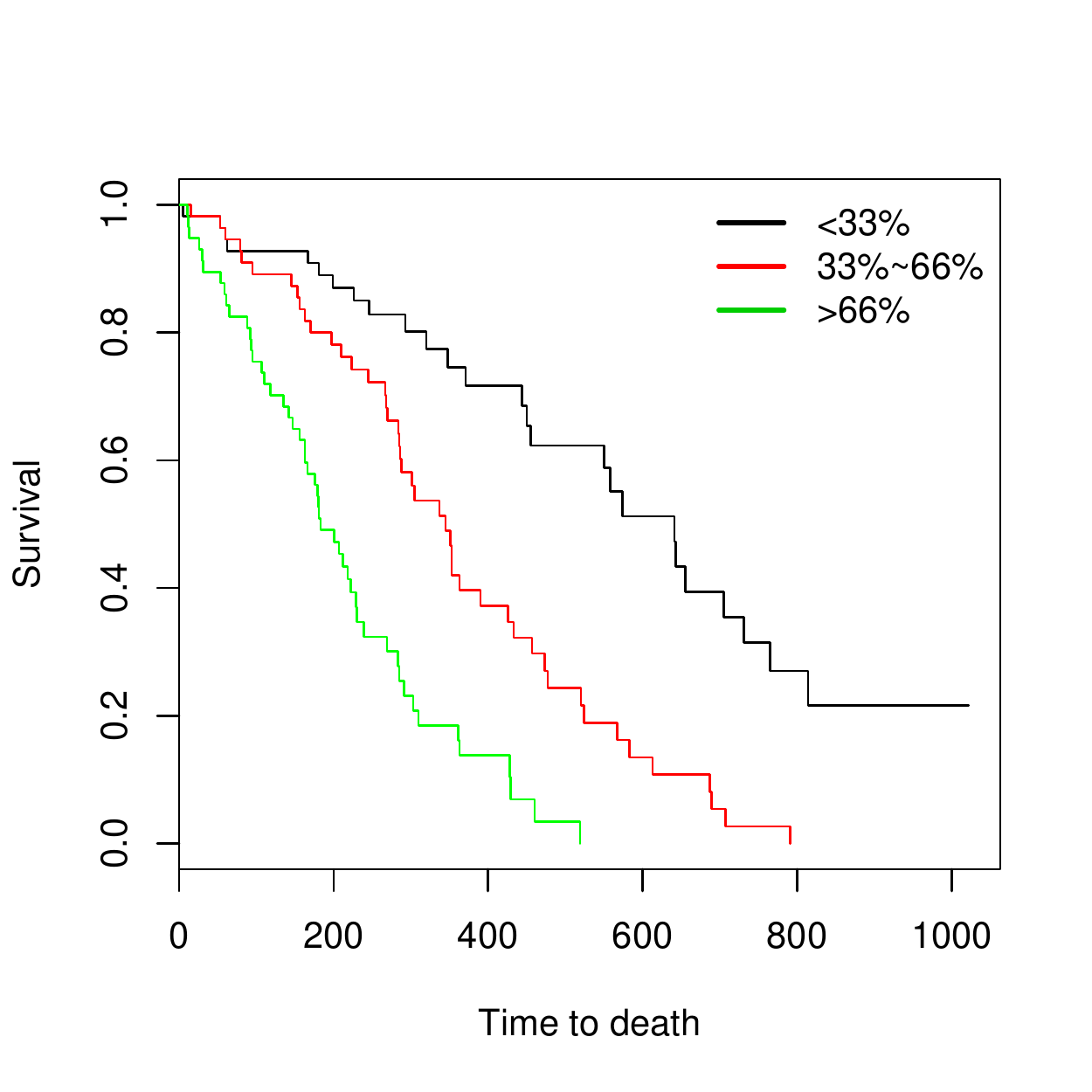}
  \caption{The Kaplan--Meier estimate of survival curves for the three risk groups  determined by the
risk scores (divided by the 33\% and 66\% quantiles) estimated from the Cox PH model with all
covariates}\label{lung.all}
\end{figure}

We fit the Cox  proportional hazards (PH)  model using these covariates and all the interaction terms. The Harrell's C-index is 0.725 which means that this model fits the data well.
We also plot the Kaplan--Meier estimate of the
survival curves for the three risk groups determined by the risk scores (divided by the 33\% and 66\% quantiles)
estimated from the Cox PH model with all covariates and interaction terms. The result is displayed in Figure \ref{lung.all}. From this figure, it can be see that the three risk curves are well separated, which also suggests that the Cox PH model fits the data reasonably well.

\begin{figure}[H]
\begin{minipage}{0.5\linewidth}
\centering
\includegraphics[scale=0.6]{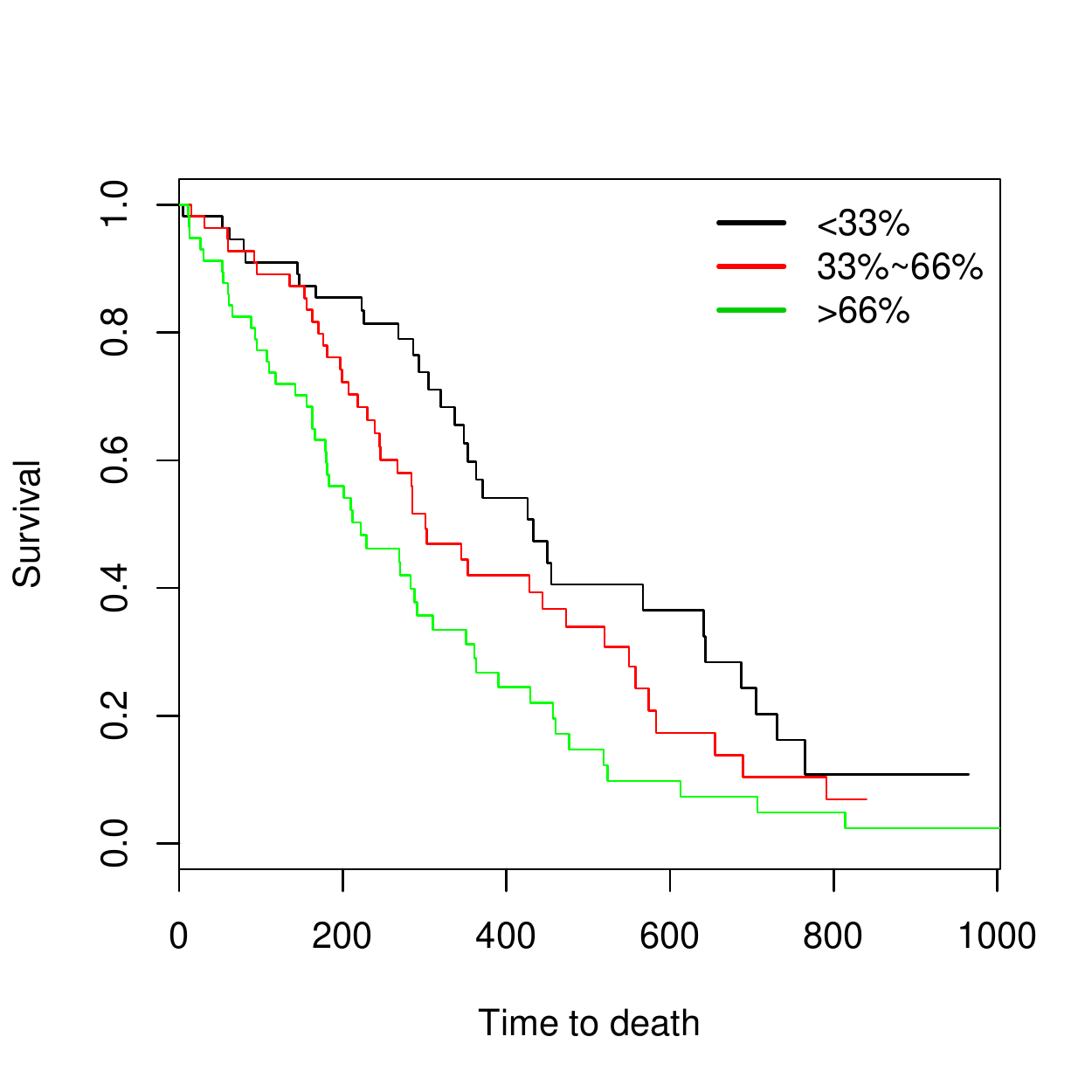}
\end{minipage}
\begin{minipage}{0.5\linewidth}
\centering
\includegraphics[scale=0.6]{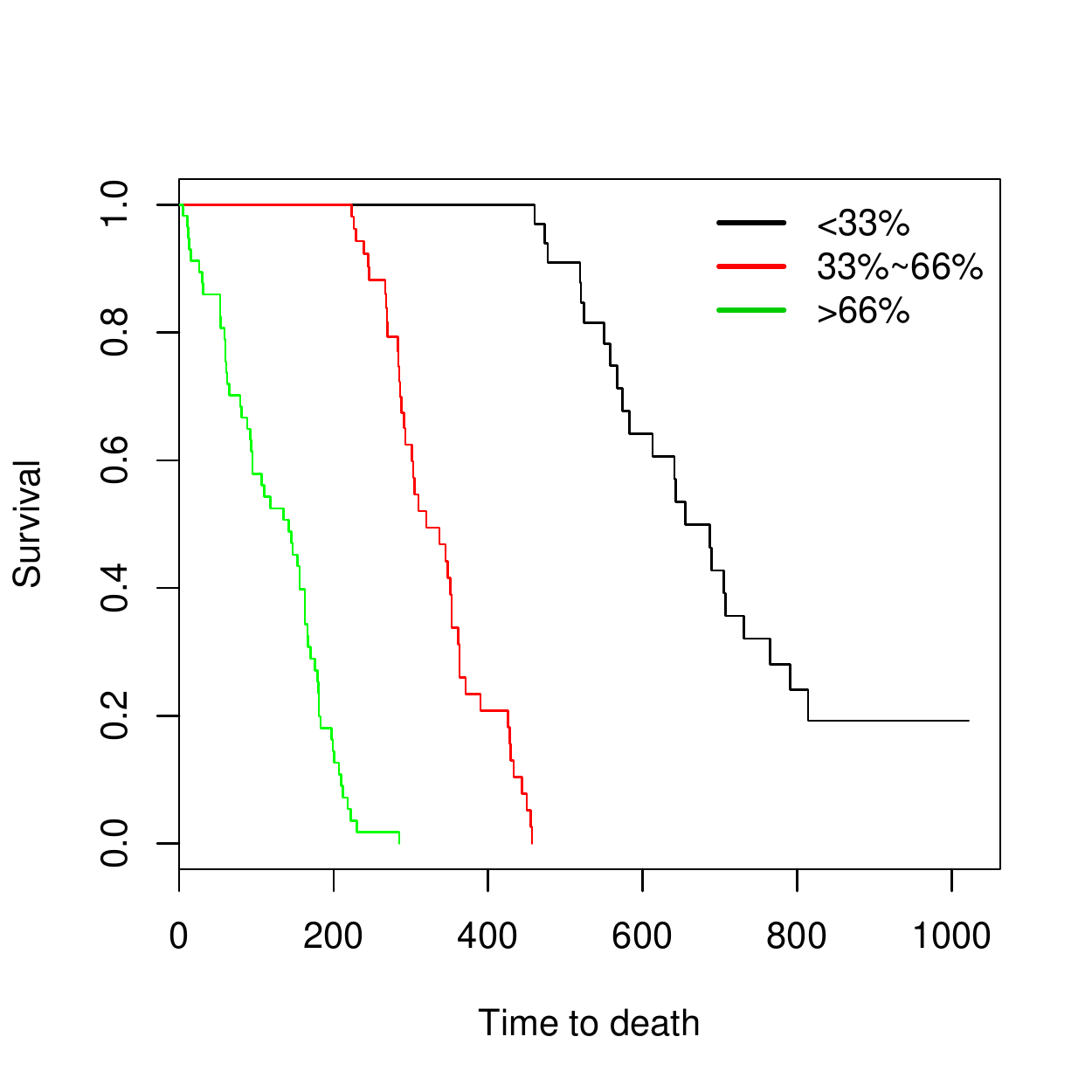}
\end{minipage}
\caption{The Kaplan--Meier estimate of survival curves for the three risk groups  determined by the
risk scores (divided by the 33\% and 66\% quantiles) estimated from the Cox PH model with DSIR directions (left) and
RDSIR directions (right).}
\label{lung.directions}
\end{figure}

For applying the proposed method to the data, we choose the Gaussian kernel, and set the regularity parameter $\tau=1$. We denote $\widehat{\nu}_1,\widehat{\nu}_2$ as the first two directions or components found by DSIR or RDSIR. The Cox PH model we fit is $\lambda(t)=\lambda_0(t)\exp\{\theta_1\nu_1+\theta_2\nu_2+\theta_3\nu_1\nu_2\}$, where $\lambda(t)$ is hazard function and $\lambda_0(t)$ is baseline hazard function. The results are displayed in Figure \ref{lung.directions}.
From this figure, we can see that the three risk group curves corresponding to RDSIR directions
show a better separation than the ones corresponding to DSIR directions, which indicates that RDSIR outperforms DSIR.

{
\section{Conclusion}
In this paper, we give a new nonlinear SDR method for survival data and demonstrate its superiority over DSIR via numerical studies. Our method is based on DSIR \cite{Li1999} and RKHS theory. We adopt the double slicing procedure to estimate weight function which is an important step for modeling survival data. The estimation of weight function requires nonparametric smoothing. We tackle it by using kernel smoothing in this paper for the convenience of theoretical development. However, other nonparametric smoothing methods such as spline \cite{spline,Wahba1990} and wavelets \cite{wavelets2009,wavelets2012} can also be used.
Another idea for handling survival data is to impute the censored $T$ observations first and then apply the SDR methods directly to the imputed data.

The choice of reproducing kernel is also worth taking into account. In general settings, we are not limited to use a single kernel. Alternatively, we can learn from the idea of multiple kernel learning  \cite{multikernel2006,multikernel2011}. As an example, consider a linear combination of $m$ kernels, that is, $K'=\sum_{i=1}^{m}K_i$ and this new function is still a kernel due to the properties of RKHS. Relevant theories and experiments of this topic deserve further study.

Other inverse regression based methods such as SAVE \cite{COOK1991} can also be generalized to nonlinear SDR for  survival data in a similar way. Other asymptotic properties such as asymptotic normality need also be further investigated. The case
where $p$ diverges with $n$ at some rate will be considered in our future research.
}
\newpage
\section{Appendix}\label{sec:Proof}

\noindent{\bf A. Proof of Lemma 2}

\renewcommand{\theequation}{A.\arabic{equation}}
\setcounter{equation}{0}

Denote $\mathcal{H}_{\Phi}=\text{span}\{\Phi\}$ and let $\bbeta$ be the eigenfunction which satisfies \eqref{eq7}.
We have shown that in \eqref{nnn}-\eqref{nnnn},
\be\widehat{\Sigma}=\dfrac{1}{n}\Phi\Phi^T\quad,\quad \widehat{\Gamma}=\dfrac{1}{n}\Phi{{Q}}\Phi^T.\label{A.1}\ee
We decompose $\bbeta$ into $\bbeta=\bbeta_1+\bbeta_2$, where $\bbeta_1\in \mathcal{H}_{\Phi}$ and
$\bbeta_2\in \mathcal{H}_{\Phi}^{\perp}$. Then, we can see that
$$\widehat{\Sigma}\bbeta_2={\bf 0}\quad,\quad \widehat{\Gamma}\bbeta_2={\bf 0}.$$
While the SDR central subspace is the intersection of all SDR subspaces, thus, the eigenfunction corresponding
to a non-zero eigenvalue is of the following form:
$$\widehat{\bm{\beta}}_{j}=\sum_{i=1}^n\widehat{\bm{\alpha}}_{j}^{(i)}\bphi(\bX_i)=\Phi\widehat{\bm{\alpha}}_{j},\quad j=1,\cdots,q.$$

Assume that the number of non-zero eigenvalues of reproducing kernel $R$ is $d_k$ and $d_k\leq \infty$, then
$$\bphi(\bs)=\Big(\sqrt{a_1}\phi_1(\bs),\sqrt{a_2}\phi_2(\bs),\cdots,\sqrt{a_{d_k}}\phi_{d_k}(\bs)\Big)^T.$$

(\Rmnum{1}) If $d_k< \infty$, $\mathcal{H}$ is finite dimensional and $\Phi$ is a $d_k\times n$ matrix.
Write the singular value decomposition (SVD) of $\Phi$ as
\begin{align}\label{A.2}
  \Phi & = UDV \notag\\
       & = (u_1,\cdots,u_{d_k})\left[
                                 \begin{array}{cc}
                                   \widebar{D}_{d\times d} &  0_{d\times (n-d)}\\
                                   0_{(d_k-d)\times d} & 0_{(n-d_k)\times(n-d)} \\
                                 \end{array}
                               \right]\left[
                                        \begin{array}{c}
                                          v_1^T \\
                                           \vdots\\
                                          v_n^T \\
                                        \end{array}
                                      \right]\\
       & = \widebar{U}\widebar{D}\widebar{V}^T,\notag
\end{align}
where $\widebar{U}=(u_1,\cdots,u_d),\ \widebar{V}=(v_1,\cdots,v_d),\ \widebar{D}=\widebar{D}_{d\times d}$ is a diagonal matrix and
$\widebar{U}^T\widebar{U}=\widebar{V}^T\widebar{V}=I_d$.

$(\Longrightarrow)$ If $\widehat{\Gamma}\bbeta=\lambda\widehat{\Sigma}\bbeta$, we have
$\big\langle \bphi(\bm{X}_i),\widehat{\Gamma}\bbeta\big\rangle=\big\langle \bphi(\bm{X}_i),\lambda\widehat{\Sigma}\bbeta\big\rangle$ for $i=1,\cdots,n$. Then
$$\Phi^T\left(\dfrac{1}{n}\Phi{{Q}}\Phi^T\bbeta\right)=\dfrac{1}{n}\Phi^T\Phi{{Q}}\Phi^T\Phi\balpha
=\dfrac{1}{n}{{R}}{{Q}}{{R}}\balpha,$$
$$\Phi^T\left(\dfrac{1}{n}\lambda\Phi\Phi^T\bbeta\right)=\dfrac{1}{n}\lambda\Phi^T\Phi\Phi^T\Phi\balpha
=\lambda\dfrac{1}{n}{{R}}^2\balpha,$$
and hence ${{R}}{{Q}}\bm{R\alpha}=\lambda{{R}}^2\bm{\alpha}$.

$(\Longleftarrow)$ If ${{R}}{{Q}}\bm{R\alpha}=\lambda{{R}}^2\bm{\alpha}$, $\bbeta=\Phi\balpha$ and
${{R}}=\Phi^T\Phi=\widebar{V}\widebar{D}^2\widebar{V}^T$, then
\be\widebar{V}\widebar{D}^2\widebar{V}^T{{Q}}\Phi^T\Phi\balpha=\lambda\widebar{V}\widebar{D}^2\widebar{V}^T
\Phi^T\Phi\balpha.\label{A.3}\ee
Multiply both sides of the equation \eqref{A.3} by $\widebar{U}\widebar{D}^{-1}\widebar{V}^T$, we obtain that
$$\widebar{U}\widebar{D}\widebar{V}^T{{Q}}\Phi^T\bbeta=\lambda\widebar{U}\widebar{D}\widebar{V}^T
\Phi^T\bbeta.$$
By \eqref{A.1} and \eqref{A.2}, we can derive that $\widehat{\Gamma}\bbeta=\lambda\widehat{\Sigma}\bbeta$.

(\Rmnum{2}) If $d_k=\infty$, then $\mathcal{H}$ is infinite dimensional. We define $\Phi$, $\Phi^T$ and
 SVD of $\Phi$ as follows: $\Phi$ is an operator from $\mathbb{R}^n$ to $\mathcal{H}$ defined
 by $\Phi\bm{\gamma}=\sum\limits_{i=1}^n\gamma_i\phi(\bm{X}_i)$ for $\bm{\gamma}=(\gamma_1,\cdots,\gamma_n)^T\in \mathbb{R}^n$ and $\Phi^T$ is the self-adjoint operator of $\Phi$. Thus, $\Phi^T$ is an operator from $\mathcal{H}$
 to $\mathscr{R}^n$ such that
 $$\Phi^Tf=\Big(\left\langle\phi(\bm{X}_1),f\right\rangle,\cdots,\left\langle\phi(\bm{X}_n),f\right\rangle\Big)^T,
 \quad\quad f\in \mathcal{H}.$$
Since the rank of compact operator $\widehat{\Sigma}$ is less than $n$, $\widehat{\Sigma}$ has the following
representation:
\be\widehat{\Sigma}=\sum_{i=1}^{d_k}\sigma_iu_i\otimes u_i=\sum_{i=1}^{d}\sigma_iu_i\otimes u_i,\label{A.4}\ee
where $d\leq n$ and $\sigma_1\geq \cdots\geq \sigma_d\geq \sigma_{d+1}=\cdots=0$. Denote $\widebar{U}=(u_1,\cdots,u_d)$.
Similar to $\Phi$, we can formally define operators $\widebar{U}$ and $\widebar{U}^T$.

By \eqref{A.4}, we know that $\phi(\bm{X}_i)\in \text{span}\{u_1,\cdots,u_d\}$. Thus, we denote $\phi(\bm{X}_i)=\widebar{U}\omega_i,
\ \Omega=(\omega_1,\cdots,\omega_n)^T$. It is easy to check that $\Omega^T\Omega=\text{diag}(n\sigma_1,\cdots,\sigma_d)$.
Denote
$$\widebar{D}_{d\times d}=\text{diag}(\sqrt{n\sigma_1},\cdots,\sqrt{n\sigma_d})\quad,\quad \widebar{V}=\Omega\widebar{D}^{-1}.$$
Then we can define SVD of $\Phi$:
$$\Phi=\widebar{U}\widebar{D}\widebar{V}^T.$$
Subsequent derivation is similar to the case that $d_k< \infty$. $\hfill \square$

\newpage
\noindent{\bf B. Proof of Lemma 3}
\renewcommand{\theequation}{B.\arabic{equation}}
\setcounter{equation}{0}

We first prove that
\be S_T(t|\bphi(\bm{X}))=\exp\left\{-E\left[\dfrac{I(t>\widetilde{T},\Delta=1)}{S_{\widetilde{T}}
(\widetilde{T}|\bphi(\bm{X}))}\bigg|\bphi(\bm{X})\right]\right\}.\label{B.1}\ee
By the condition $1$, $T\independent C\mid \bX$, it follows that
\begin{equation*}S_{\widetilde{T}}(t|\bphi(\bm{X}))=S_T(t|\bphi(\bm{X}))S_C(t|\bphi(\bm{X})),\end{equation*}
where $S_C(t|\bphi(\bm{X}))=P\{C>t|\bphi(\bm{X})\}$. Thus,
\begin{align*}
  E\left\{\dfrac{I(t>\widetilde{T},\Delta=1)}{S_{\widetilde{T}}
  (\widetilde{T}|\bphi(\bm{X}))}\bigg|\bphi(\bm{X})\right\} & =E\left\{\dfrac{I(t>T)I(T<C)}{S_T
(T|\bphi(\bm{X}))S_C(T|\bphi(\bm{X}))}\bigg|\bphi(\bm{X})\right\} \nonumber\\
   & =E\left\{\dfrac{I(t>T)}{S_T
(T|\bphi(\bm{X}))S_C(T|\bphi(\bm{X}))}E(I(T<C)|\bphi(\bm{X}),T)\bigg|\bphi(\bm{X})\right\} \nonumber\\
   & =E\left\{\dfrac{I(t>T)}{S_T(T|\bphi(\bm{X}))}\bigg|\bphi(\bm{X})\right\}.
\end{align*}
Let $\lambda(s)=f(s)/S_T(s)$, the hazard function. Since the cumulative hazard function $H(t)$ satisfies that
\be\nonumber H(t)=\int_{0}^{t}\lambda(s)ds=\int_{0}^{t}f(s)/S_T(s)ds=E\left[\dfrac{I(T<t)}{S_T(T)}\right],\ee
by the fact that $S_T(t)=\exp(-H(t))$, \eqref{B.1} follows.
By the direct calculation, we obtain that
\begin{align*}
  w(t',t,\bphi(\bm{X})) & =\dfrac{S_T(t|\bphi(\bm{X}))}{S_T(t'|\bphi(\bm{X}))} \no\\
    & = \exp\left\{\dfrac{I(t>\widetilde{T},\Delta=1)-I(t'>\widetilde{T},\Delta=1)}{S_{\widetilde{T}}
(\widetilde{T}|\bphi(\bm{X}))}\bigg|\bphi(\bm{X})\right\}\no\\
   & = \exp\left\{\dfrac{I(t'\leq\widetilde{T}<t,\Delta=1)}{S_{\widetilde{T}}
(\widetilde{T}|\bphi(\bm{X}))}\bigg|\bphi(\bm{X})\right\}\no\\
   & = \exp\left\{-\Lambda(t',t|\bphi(\bm{X}))\right\}.
\end{align*}
Since $\widetilde{T}=T\wedge C$, it follows that $\{\widetilde{T}\geq t\}=\{T\geq t, C\geq t\}$. Hence, we have
\begin{align*}
  E\{\bphi(\bm{X})I(T\geq t)\} & = E\{\bphi(\bm{X})I(T\geq t, C\geq t)\}+E\{\bphi(\bm{X})I(T\geq t, C< t)\}\no\\
   & = E\{\bphi(\bm{X})I(\widetilde{T} \geq t)\}+E\{\bphi(\bm{X})I(T\geq t, C< t)\}. \label{(B.1)}
\end{align*}
Further, we obtain that
\begin{align*}
  E\{\bphi(\bm{X})I(T\geq t, C< t)\} & = E\{\bphi(\bm{X})I(\widetilde{T}< t, \Delta=0)I(T\geq t)\}\no\\
   & = E\{\bphi(\bm{X})I(T\geq t, C< t)E[I(T\geq t)|\widetilde{T},\Delta=0,\bphi(\bm{X})]\}\no\\
   & = E\{\bphi(\bm{X})I(T\geq t, C< t)E[I(T\geq t)|C,T>C,\bphi(\bm{X})]\}\no\\
   & = E\{\bphi(\bm{X})I(\widetilde{T}< t, \Delta=0)w(C,t,\bphi(\bm{X}))\}\no\\
   & = E\{\bphi(\bm{X})I(\widetilde{T}< t, \Delta=0)w(\widetilde{T},t,\bphi(\bm{X}))\}.
\end{align*}

The derivation for $E\{I(T>t)\}$ is similar to the above. Thus, Lemma \ref{conditionsurv} is proved.$\hfill{} \square$

\noindent{\bf C. Proof of Lemma 4}

\renewcommand{\theequation}{C.\arabic{equation}}
\setcounter{equation}{0}

Before we prove Lemma \ref{lemma4}, we first present four propositions.

\renewcommand{\theproposition}{C.\arabic{proposition}}

\setcounter{proposition}{0}

Denote
\be\widetilde{P}(T\geq t)=\dfrac{1}{n}\sum\limits_{i=1}^{n}I(\widetilde{T}_i\geq t)+\dfrac{1}{n}\sum\limits_{i:\widetilde{T}_i<t,\triangle_i=0}\widehat{w}\big(\widetilde{T}_i,t,\bm{Z}_i\big).\label{C.1}\ee

\begin{proposition}\label{prop1} Under {Assumptions \eqref{A1}-\eqref{A5}}, we have
$\Big|\widetilde{P}(T\geq t)-P(T\geq t)\Big|=O_p(n^{-1/2}).$

\end{proposition}

{\bf Proof}. Denote $I_i=I(\widetilde{T}_i<t,\triangle_i=0)$. Thus, by \eqref{www}, the second term of the right hand side of \eqref{C.1} has the following expression:
\begin{equation}
\dfrac{1}{n}\sum\limits_{i:\widetilde{T}_i<t,\triangle_i=0}\widehat{w}\big(\widetilde{T}_i,t,\bm{Z}_i\big)=\dfrac{1}{n}\sum\limits_iI_i\ \widehat{w}(\widetilde{T}_i,t,\bm{Z}_i)=\dfrac{1}{n}\sum\limits_iI_i\ e^{-\widehat{\Lambda}(\widetilde{T}_i,t|\bm{Z}_i)}.\label{C.19n}
\end{equation}
By \eqref{llll}, we have
\begin{equation}\widehat{\Lambda}_i=\widehat{f}(\bm{Z}_i)^{-1}n^{-1}\sum_jI_{ij}
\widehat{S}_{\widetilde{T}}(\widetilde{T}_j|\bm{Z}_j)^{-1}e_{ji}.\label{yyy}\end{equation}
In terms of Assumptions \eqref{A2}-\eqref{A4}, by using the Taylor's expansion, we have
\begin{align}
  \widehat{S}_{\widetilde{T}}(\widetilde{T}_j|\bm{Z}_j)^{-1} &
  =\widehat{f}(\bm{Z}_j)\Big[n^{-1}\sum_kI(\widetilde{T}_k>\widetilde{T}_j)e_{kj}\Big]^{-1}\no\\
  & = \widehat{f}(\bm{Z}_j)\left\{\dfrac{1}{n}\sum_kv_{kj}+
  \dfrac{1}{n}\sum_kE\Big[I(\widetilde{T}_k>\widetilde{T}_j)e_{kj}\big|\bm{Z}_j,\widetilde{T}_j\Big]\right\}^{-1}\no\\
  & = \widehat{f}(\bm{Z}_j)\Big[S_{\widetilde{T}}(\widetilde{T}_j|\bm{Z}_j)f(\bm{Z}_j)
  +\dfrac{1}{n}\sum_kv_{kj}+O_p(n^{-1/2})\Big]^{-1} \no\\
  & = \widehat{f}(\bm{Z}_j)\Big[f(\bm{Z}_j)^{-1}S_{\widetilde{T}}(\widetilde{T}_j|\bm{Z}_j)^{-1}
  -f(\bm{Z}_j)^{-2}S_{\widetilde{T}}(\widetilde{T}_j|\bm{Z}_j)^{-2}\dfrac{1}{n}\sum_kv_{kj}+O_p(n^{-1/2})\Big]\no\\
    & =S_{\widetilde{T}}(\widetilde{T}_j|\bm{Z}_j)^{-1}-f(\bm{Z}_j)^{-1}
  S_{\widetilde{T}}(\widetilde{T}_j|\bm{Z}_j)^{-2}\dfrac{1}{n}\sum_kv_{kj}\no\\
  & \quad + f(\bm{Z}_j)^{-1}S_{\widetilde{T}}(\widetilde{T}_j|\bm{Z}_j)^{-1}\dfrac{1}{n}\sum_ku_{kj}+O_p(n^{-1/2}).\label{C.7}
\end{align}
By Assumption \eqref{A1}, it follows that
\be\label{C.3}\widehat{f}(\bm{Z}_i)=f(\bm{Z}_i)+\dfrac{1}{n}\sum_{k}u_{ki}+O_p(n^{-1/2}),\ee
which implies that
\begin{align}
  \widehat{f}(\bm{Z}_i)^{-1} & = f_i^{-1}-f_i^{-2}\left[\dfrac{1}{n}\sum_ku_{ki}+O_p(n^{-1/2})\right]+O_p(n^{-1/2})\no\\
     & = f_i^{-1}-f_i^{-2}\dfrac{1}{n}\sum_ku_{ki}+O_p(n^{-1/2}).\label{C.8}
\end{align}
Combining \eqref{yyy}, \eqref{C.7} and \eqref{C.8},  we have
\begin{align}
  \widehat{\Lambda}_i & = \widehat{f}(\bm{Z}_i)^{-1}n^{-1}\sum_jI_{ij}
  e_{ji}\bigg[S_j^{-1}-f_j^{-1}S_j^{-2}n^{-1}\sum_kv_{kj}\no\\
  & \qquad\qquad+ f_j^{-1}S_j^{-1}n^{-1}\sum_ku_{kj}+O_p(n^{-1/2})\bigg]\no\\
  & =\widehat{f}(\bm{Z}_i)^{-1}n^{-1}\sum_jI_{ij}e_{ji}S_j^{-1}-\widehat{f}(\bm{Z}_i)^{-1}n^{-1}\sum_jI_{ij}e_{ji}
  f_j^{-1}S_j^{-2}n^{-1}\sum_kv_{kj}\no\\
  & \qquad\qquad+\widehat{f}(\bm{Z}_i)^{-1}n^{-1}\sum_jI_{ij}e_{ji}f_j^{-1}S_j^{-1}n^{-1}\sum_ku_{kj}\no\\
  &\qquad\qquad+\widehat{f}(\bm{Z}_i)^{-1}n^{-1}\sum_jI_{ij}e_{ji}O_p(n^{-1/2})\no\\
  &={\tt Q}_{11}-{\tt Q}_{12}+{\tt Q}_{13}+{\tt Q}_{14},\label{C.9}
\end{align}
and we consider each of these four terms in the following.

{By \eqref{C.8}, we have
\begin{align}
  {\tt Q}_{11} & = \widehat{f}(\bm{Z}_i)^{-1}n^{-1}\sum_jI_{ij}e_{ji}S_j^{-1}\no\\
              & = \Big[f_i^{-1}-f_i^{-2}\dfrac{1}{n}\sum_ku_{ki}+O_p(n^{-1/2})\Big]n^{-1}\sum_jI_{ij}e_{ji}S_j^{-1}\no\\
         & = f_i^{-1}n^{-1}\sum_jI_{ij}e_{ji}S_j^{-1}-f_i^{-2}n^{-1}\sum_jI_{ij}e_{ji}S_j^{-1}
         \dfrac{1}{n}\sum_ku_{ki}+n^{-1}\sum_jI_{ij}e_{ji}S_j^{-1}O_p(n^{-1/2})\no\\
         &={\tt Q}_{21}-{\tt Q}_{24}+O_p(n^{-1/2}).\label{C.10}
\end{align}
The last term $O_p(n^{1/2})$ in \eqref{C.10} is obtained by the application of the law of large numbers on
 $n^{-1}\sum_jI_{ij}e_{ji}S_j^{-1}O_p(n^{-1/2})$.

Similarly, we have
\begin{align}
  {\tt Q}_{12} & = \widehat{f}(\bm{Z}_i)^{-1}n^{-1}\sum_jI_{ij}e_{ji}f_j^{-1}S_j^{-2}\dfrac{1}{n}\sum_kv_{kj}\no\\
  &=\Big[f_i^{-1}-f_i^{-2}\dfrac{1}{n}\sum_ku_{ki}+O_p(n^{-1/2})\Big]n^{-1}\sum_jI_{ij}e_{ji}f_j^{-1}S_j^{-2}\dfrac{1}{n}\sum_kv_{kj}\no\\
&=f_i^{-1}n^{-1}\sum_jI_{ij}e_{ji}f_j^{-1}S_j^{-2}\dfrac{1}{n}\sum_kv_{kj}-f_i^{-2}n^{-1}\sum_jI_{ij}e_{ji}f_j^{-1}S_j^{-2}\Big(\dfrac{1}{n}\sum_kv_{kj}\Big)\Big(\dfrac{1}{n}\sum_ku_{kj}\Big)\no\\
&\qquad+n^{-1}\sum_jI_{ij}e_{ji}f_j^{-1}S_j^{-2}\dfrac{1}{n}\sum_kv_{kj}O_p(n^{-1/2})\no\\
& = f_i^{-1}n^{-1}\sum_jI_{ij}e_{ji}f_j^{-1}S_j^{-2}\dfrac{1}{n}\sum_kv_{kj}+O_p(n^{-1/2})\no\\
&={\tt Q}_{22}+O_p(n^{-1/2}).\label{C.11}
\end{align}
The fourth ``$=$" is obtained as follows: The applications of Assumptions (A2) and and (A4) yield that $\big(\frac{1}{n}\sum_kv_{kj}\big)\big(\frac{1}{n}\sum_ku_{kj}\big)=O_p(n^{-1/2})$, and the application of the law of large numbers gives the subsequent conclusion.

We can derive ${\tt Q}_{13}$ and ${\tt Q_{14}}$ in the same way as for ${\tt Q}_{11}$ and ${\tt Q_{12}}$, and obtain that
\begin{align}
  {\tt Q}_{13} & = \widehat{f}(\bm{Z}_i)^{-1}n^{-1}\sum_jI_{ij}e_{ji}f_j^{-1}S_j^{-1}\dfrac{1}{n}\sum_ku_{kj}\no\\
         & = f_i^{-1}n^{-1}\sum_jI_{ij}e_{ji}f_j^{-1}S_j^{-1}\dfrac{1}{n}\sum_ku_{kj}+O_p(n^{-1/2})\no\\
         &={\tt Q}_{23}+O_p(n^{-1/2}),\label{C.12}
\end{align}
and
\be {\tt Q}_{14}=\widehat{f}(\bm{Z}_i)^{-1}n^{-1}\sum_jI_{ij}e_{ji}O_p(n^{-1/2})=O_p(n^{-1/2}).\label{C.13}\ee}
Combining \eqref{C.9}$-$\eqref{C.13}, we obtain that
\be\widehat{\Lambda}_i={\tt Q}_{21}-{\tt Q}_{22}+{\tt Q}_{23}-{\tt Q}_{24}+O_p(n^{-1/2}),\label{C.14}\ee
where
\begin{align*}
  {\tt Q}_{21} & = \dfrac{1}{nf_i}\sum_jI_{ij}e_{ji}S_j^{-1},\\
  {\tt Q}_{22} & = \dfrac{1}{nf_i}\sum_jI_{ij}e_{ji}f_j^{-1}S_j^{-2}\dfrac{1}{n}\sum_kv_{kj},\\
  {\tt Q}_{23} & = \dfrac{1}{nf_i}\sum_jI_{ij}e_{ji}f_j^{-1}S_j^{-1}\dfrac{1}{n}\sum_ku_{kj},\\
  {\tt Q}_{24} & = \dfrac{1}{nf_i^2}\sum_jI_{ij}e_{ji}S_j^{-1}\dfrac{1}{n}\sum_ku_{ki}.
\end{align*}
Thus, by \eqref{C.19n}%, \eqref{C.12}
and \eqref{C.14}, we have
\begin{align}
&\dfrac{1}{n}\sum\limits_iI_i\ \widehat{w}(\widetilde{T}_i,t,\bm{Z}_i)=\dfrac{1}{n}\sum\limits_iI_i\ e^{-\widehat{\Lambda}(\widetilde{T}_i,t|\bm{Z}_i)}\nonumber\\
=&\dfrac{1}{n}\sum\limits_iI_i
e^{-\left({\tt Q}_{21}-{\tt Q}_{22}+{\tt Q}_{23}-{\tt Q}_{24}+O_p(n^{-1/2})\right)}\label{dldl}
\end{align}
To examine ${\tt Q}_{21}$, we define \be\varepsilon_{ij}=I_{ij}S_j^{-1}e_{ji}-E\left(I_{ij}S_j^{-1}e_{ji}|\bm{Z}_i,\widetilde{T}_i\right).\label{C.16}\ee
{By  Li \textit{et al.} \cite{Li1999}, $\varepsilon_{ij}$ satisfies that
\be\label{eeij}
E\varepsilon_{ij}=0, \quad E(\varepsilon_{ij}|\bm{Z}_i,\widetilde{T}_i)=0,\quad\text{var}(\varepsilon_{ij})=O(h_n^{-m}),
\ee}
In view of \eqref{C.16} and Assumption \eqref{A5}, we have
\begin{align}
{\tt Q}_{21}&=\dfrac{1}{f_i}\left[\dfrac{1}{n}\sum_j\varepsilon_{ij}+
\dfrac{1}{n}\sum_jE\left(I_{ij}S_j^{-1}e_{ji}|\bm{Z}_i,\widetilde{T}_i\right)\right]\nonumber\\
& =\dfrac{1}{f_i}\left[\dfrac{1}{n}\sum_j\varepsilon_{ij}+\Lambda_if_i+O_p(n^{-1/2})\right]\nonumber\\
& = \Lambda_i+\dfrac{1}{nf_i}\sum_j\varepsilon_{ij}+O_p(n^{-1/2}),\label{C.17}
\end{align}
which, jointly with \eqref{dldl}, yields that
\begin{align}
&\dfrac{1}{n}\sum\limits_iI_i\ \widehat{w}(\widetilde{T}_i,t,\bm{Z}_i)=\dfrac{1}{n}\sum\limits_iI_ie^{-\Lambda_i}
e^{-\left(f_i^{-1}n^{-1}\sum\limits_j\varepsilon_{ij}-{\tt Q}_{22}+{\tt Q}_{23}-{\tt Q}_{24}+O_p(n^{-1/2})\right)}\nonumber\\
=&\dfrac{1}{n}\sum\limits_iI_ie^{-\Lambda_i}
\left[1-\left(\dfrac{1}{nf_i}\sum\limits_j\varepsilon_{ij}-{\tt Q}_{22}+{\tt Q}_{23}-{\tt Q}_{24}+O_p(n^{-1/2})\right)(1+o_p(1))\right]\nonumber\\
=&\dfrac{1}{n}\sum\limits_iI_iw_i-\left\{\dfrac{1}{n^2}\sum\limits_{i,j}I_iw_if_i^{-1}\varepsilon_{ij}-\dfrac{1}{n}\sum\limits_iI_iw_i
\Big[{\tt Q}_{22}-{\tt Q}_{23}+{\tt Q}_{24}+O_p(n^{-1/2})\Big]\right\}(1+o_p(1))\nonumber\\
=&{\tt Q}_{31}-({\tt Q}_{32}-{\tt Q}_{33})(1+o_p(1)).\label{C.19}
\end{align}
We need to evaluate ${\tt Q}_{31}$. {As $I_i=I(\widetilde{T}_i<t,\triangle_i=0)$, $\Lambda_i=\Lambda(\widetilde{T}_i,t|\mathbf{P}_{\mathcal{B}_J}\bphi(\bX_i))$ and $w_i=e^{-\Lambda_i}$,  it follows that $I_iw_i$ are i.i.d., which, jointly with the central limit theorem, yields that
\be\frac{\sum\limits_iI_iw_i-E\big(\sum\limits_iI_iw_i\big)}{\sqrt{nVar(I_iw_i)}}\stackrel{d}{\longrightarrow}N(0,1).\nonumber\ee
Therefore,
\be {\tt Q}_{31}-E{\tt Q}_{31}=\dfrac{1}{n}\sum\limits_i\Big[I_iw_i-E(I_iw_i)\Big]=O_p(n^{-1/2}).\label{C.20}\ee}
Since $\Delta_i=I(T_i\leq C)$ are i.i.d. and $\widetilde{T}_i=T_i\wedge C$ are also i.i.d., by  direct calculation,  we obtain that
\begin{align*}
  E(I_iw_i) & =E\left[I(\widetilde{T}<t,\triangle=0)e^{-\Lambda(\widetilde{T},t|\bm{Z})}\right] \\
   & =E\left[I(\widetilde{T}<t,\triangle=0)w(\widetilde{T},t,\bm{Z})\right]\\
   & =E\left[I(\widetilde{T}<t,\triangle=0)w(C,t,\bm{Z})\right],
\end{align*}
where $\bm{Z}=\mathbf{P}_{\mathcal{B}_J}\bphi(\bX)$ defined in \eqref{notations}.
We can also derive that
\begin{align*}
w(C,t,\bm{Z})&=\dfrac{P(T\geq t|C,\bm{Z})}{P(T\geq C|C,\bm{Z})}
=\dfrac{P(T\geq t,T\geq C|C,\bm{Z})}{P(T\geq C|C,\bm{Z})}+
\dfrac{P(T\geq t,T < C|C,\bm{Z})}{P(T\geq C|C,\bm{Z})}\\
&=P(T\geq t|C,T\geq C,\bm{Z}),
\end{align*}
and hence we have
\begin{align}
  E(I_iw
  _i) & = E\left[I(\widetilde{T}<t,\triangle=0)P(T\geq t|C,T\geq C,\bm{Z})\right]\no\\
   & =E\left[I(\widetilde{T}<t,\triangle=0)P(T\geq t|\widetilde{T},\triangle=0,\bm{Z})\right]\no\\
   & =E\left[I(\widetilde{T}<t,\triangle=0)E\left(I(T\geq t)\Big|\widetilde{T},\triangle=0,\bm{Z}\right)\right]\no\\
   & =E\left\{E\left[I(\widetilde{T}<t,\triangle=0)I(T\geq t)\Big|\widetilde{T},\triangle=0,\bm{Z}\right]\right\}\no\\
   & =E\left[I(\widetilde{T}<t,\triangle=0)I(T\geq t)\right]\no\\
   & =P\left(\widetilde{T}<t,\triangle=0,T\geq t\right)\no\\
   & =P\left(T\geq t,C<t\right). \label{C.21}
\end{align}
which, combined with \eqref{C.20}, yields that
\be {\tt Q}_{31}=P\left(T\geq t,C<t\right)+O_p(n^{-1/2}).\label{jkjk}\ee
Next we evaluate ${\tt Q}_{32}$. In light of \eqref{eeij}, it follows that
\[E(I_iw_if_i^{-1}\varepsilon_{ij}I_{i'}w_{i'}f_{i'}^{-1}\varepsilon_{i'j'})
=\left\{\begin{array}{ll}
          0, & \hbox{$j\neq j'$,} \\
       O(1), & \hbox{$j=j',i\neq i'$,}\\
       O(h_n^{-m}), & \hbox{$j=j',i= i'$,}
        \end{array}
  \right.
\]
therefore \be \text{var}\left(n^{-2}\sum\limits_{i,j}I_iw_if_i^{-1}\varepsilon_{ij}\right)=n^{-4}[n^3O(1)+n^2O(h_n^{-m})]=O\left(\dfrac{1}{n}\right),\label{C.22}\ee
which implies that
\be{\tt Q}_{32}=O_p(n^{-1/2}).\label{jkjk1}\ee
We now assess ${\tt Q}_{33}$.  We examine each of its three terms in sequence. We have
\begin{align*}
  n^{-1}\sum\limits_iI_iw_i{\tt Q}_{22} & = n^{-1}\sum_iI_iw_if_i^{-1}n^{-1}
\sum_jI_{ij}e_{ji}f_j^{-1}S_j^{-2}n^{-1}\sum_kv_{kj} \\
   &= n^{-2}\sum_j\left[\sum_iI_iw_if_i^{-1}I_{ij}f_j^{-1}S_j^{-2}e_{ji}\right]n^{-1}\sum_kv_{kj}\\
   & \xlongequal{\triangle} n^{-2}\sum_{j,k}\widetilde{a}_jv_{kj},
\end{align*}
where $\widetilde{a}_j=n^{-1}\sum\limits_iI_iw_if_i^{-1}I_{ij}f_j^{-1}S_j^{-2}e_{ji}$.
{In light of Li \textit{et al.} \cite{Li1999}, we have
\be \text{var}\left(n^{-2}\sum_{j,k}\widetilde{a}_jv_{kj}\right)=O(n^{-1}).\label{C.23}\ee
}
We can evaluate the variance of $n^{-1}\sum\limits_iI_iw_i{\tt Q}_{23}$ exactly in the same way and obtain that
\be
\text{var}(n^{-1}\sum\limits_iI_iw_i{\tt Q}_{23})=O(n^{-1}).\label{C.24}\ee
Similarly, we can show that the variance of $n^{-1}\sum\limits_iI_iw_i{\tt Q}_{24}$
also has the order of $n^{-1}$, i.e.,
\begin{align}
  \text{var}\Big(n^{-1}\sum\limits_iI_iw_i{\tt Q}_{24}\Big) & =O(n^{-1}). \label{C.25}
\end{align}
Combination of \eqref{C.23}-\eqref{C.25} yields that
\be{\tt Q}_{33}=O_p(n^{-1/2}). \label{jkjk2}\ee
We conclude from \eqref{C.19}, \eqref{jkjk}, \eqref{jkjk1}, and \eqref{jkjk2} that
\be\dfrac{1}{n}\sum_iI_i\widehat{w}(\widetilde{T}_i,t,\bm{Z}_i)=P(T\geq t,C<t)+O_p(n^{-1/2}).\label{ttt}\ee
The remaining task is to assess the first term of the right hand side of \eqref{C.1}. By the central limit theorem, it follows that
\begin{align}
  \dfrac{1}{n}\sum\limits_{i=1}^{n}I(\widetilde{T}_i\geq t) & =EI(\widetilde{T}\geq t)+O_p(n^{-1/2}) \no\\
   & =P(\widetilde{T}\geq t)+O_p(n^{-1/2})\no\\
   & =P(T\geq t,C\geq t) +O_p(n^{-1/2}).\label{jkjk3}
\end{align}
Combining \eqref{C.1}, \eqref{ttt}, and \eqref{jkjk3} yields
\begin{align}
  \widetilde{P}(T\geq t) & = \dfrac{1}{n}\sum\limits_{i=1}^{n}I(\widetilde{T}_i\geq t)
+ \dfrac{1}{n}\sum_iI_i\widehat{w}(\widetilde{T}_i,t,\bm{Z}_i)\no\\
       & =P(T\geq t,C\geq t)+P(T\geq t,C< t)+O_p(n^{-1/2})\no\\
       & =P(T\geq t) + O_p(n^{-1/2}), \label{C.26}
\end{align}
which completes the proof.
$\hfill{} \square$

\begin{proposition}\label{prop2}
Under {Assumptions \eqref{A6}-\eqref{A8}}, we have
$$\Big|\widehat{P}(T\geq t)-\widetilde{P}(T\geq t)\Big|=O_p\left(n^{-b}h_n^{-(2m+1)}\right).$$
\end{proposition}

{\bf Proof}. We only need to show that
\begin{align*}
  &\dfrac{1}{n}\sum\limits_{i=1}^nI_i\widehat{w}(\widetilde{T}_i,t,\bm{Z}_i)
-\dfrac{1}{n}\sum\limits_{i=1}^nI_i\widehat{w}(\widetilde{T}_i,t,\widehat{\bm{Z}}_i)   \\
   =& \dfrac{1}{n}\sum\limits_{i=1}^nI_i\left[\widehat{w}(\widetilde{T}_i,t,\bm{Z}_i)
-\widehat{w}(\widetilde{T}_i,t,\widehat{\bm{Z}}_i)\right]\\
=&O_p\left(n^{-b}h_n^{-(2m+1)}\right).
\end{align*}
Through direct calculation and the Taylor's expansion, we obtain that
\begin{align}
&\widehat{w}(t',t,\widehat{\bm{Z}})-\widehat{w}(t',t,\bm{Z})\no\\
=&e^{-\widehat{\Lambda}(t',t|\widehat{\bm{Z}})}-e^{-\widehat{\Lambda}(t',t|\bm{Z})}\no\\
=&e^{-\widehat{\Lambda}(t',t|\bm{Z})}\Big\{e^{-\left[\widehat{\Lambda}(t',t|\widehat{\bm{Z}})
-\widehat{\Lambda}(t',t|\bm{Z})\right]}-1\Big\}\no\\
=&e^{-\widehat{\Lambda}(t',t|\bm{Z})}\left\{1-\left[\widehat{\Lambda}(t',t|\widehat{\bm{Z}})
-\widehat{\Lambda}(t',t|\bm{Z})\right]\cdot\Big(1+o_p(1)\Big)-1\right\}\no\\
=&-e^{-\widehat{\Lambda}(t',t|\bm{Z})}\left[\widehat{\Lambda}(t',t|\widehat{\bm{Z}})
-\widehat{\Lambda}(t',t|\bm{Z})\right]\cdot\Big(1+o_p(1)\Big).\label{C.27}
\end{align}
Then it suffices to consider the term $\widehat{\Lambda}(t',t|\widehat{\bm{Z}})-\widehat{\Lambda}(t',t|\bm{Z})$. We have
\begin{align}
&\widehat{\Lambda}(t',t|\widehat{\bm{Z}})-\widehat{\Lambda}(t',t|\bm{Z})\no\\
&=\dfrac{n^{-1}\sum\limits_jI(t'<\widetilde{T}_j<t,\triangle_i=1)
\widehat{S}_{\widetilde{T}}(\widetilde{T}_j|\widehat{\bm{Z}})^{-1}h_n^{-m}
K_{m}\left(h_n^{-1}(\widehat{\bm{Z}}_j-\widehat{\bm{Z}})\right)}
{\widehat{f}(\widehat{\bm{Z}})}\no\\
&-\dfrac{n^{-1}\sum\limits_jI(t'<\widetilde{T}_j<t,\triangle_i=1)
\widehat{S}_{\widetilde{T}}(\widetilde{T}_j|\bm{Z})^{-1}h_n^{-m}
K_{m}\Big(h_n^{-1}(\bm{Z}_j-\bm{Z})\Big)}
{\widehat{f}(\widehat{\bm{Z}})}\no\\
&\xlongequal{\triangle} \dfrac{\text{\Rmnum{2}}_1}{\widehat{f}(\widehat{\bm{Z}})}-
\dfrac{\text{\Rmnum{2}}_2}{\widehat{f}(\widehat{\bm{Z}})}
=\dfrac{\text{\Rmnum{2}}_1-\text{\Rmnum{2}}_2}{f(\widehat{\bm{Z}})}
+\text{\Rmnum{2}}_1\left(\dfrac{1}{\widehat{f}(\widehat{\bm{Z}})}
-\dfrac{1}{\widehat{f}(\bm{Z})}\right),\label{C.28}
\end{align}
where
\[
\text{\Rmnum{2}}_1=n^{-1}\sum\limits_jI(t'<\widetilde{T}_j<t,\triangle_i=1)
\widehat{S}_{\widetilde{T}}(\widetilde{T}_j|\widehat{\bm{Z}})^{-1}h_n^{-m}
K_{m}\left(h_n^{-1}(\widehat{\bm{Z}}_j-\widehat{\bm{Z}})\right)\]
and
\begin{align}
&\text{\Rmnum{2}}_1-\text{\Rmnum{2}}_2\nonumber\\
=&\dfrac{1}{n}\sum\limits_jI(t'<\widetilde{T}_j<t,\triangle_i=1)
\Big[\widehat{S}_{\widetilde{T}}(\widetilde{T}_j|\widehat{\bm{Z}})^{-1}h_n^{-m}
K_{m}\left(h_n^{-1}(\widehat{\bm{Z}}_j-\widehat{\bm{Z}})\right)\nonumber\\
&\qquad\qquad\qquad\qquad\qquad\qquad\quad-\widehat{S}_{\widetilde{T}}(\widetilde{T}_j|\bm{Z})^{-1}h_n^{-m}
K_{m}\Big(h_n^{-1}(\bm{Z}_j-\bm{Z})\Big)\Big]\nonumber\\
\xlongequal{\triangle}&\dfrac{1}{n}\sum\limits_jI(t'<\widetilde{T}_j<t,\triangle_i=1)\text{\Rmnum{3}}_j,\label{C.29}\end{align}
with
\begin{align}
\text{\Rmnum{3}}_j&=\widehat{S}_{\widetilde{T}}(\widetilde{T}_j|\bm{Z})^{-1}h_n^{-m}
\left[K_{m}\left(h_n^{-1}(\widehat{\bm{Z}}_j-\widehat{\bm{Z}})\right)
-K_{m}\Big(h_n^{-1}(\bm{Z}_j-\bm{Z})\Big)\right]\no\\
\qquad\quad&+\left[\widehat{S}_{\widetilde{T}}(\widetilde{T}_j|\widehat{\bm{Z}})^{-1}
-\widehat{S}_{\widetilde{T}}(\widetilde{T}_j|\bm{Z})^{-1}\right]
h_n^{-m}K_{m}\left(h_n^{-1}(\widehat{\bm{Z}}_j-\widehat{\bm{Z}})\right).\label{C.30}
\end{align}
\noindent By \eqref{C.28}-\eqref{C.30}, we need to find out the convergence rates of the following terms
\begin{align}&\widehat{f}(\widehat{\bm{Z}})^{-1}-\widehat{f}(\bm{Z})^{-1},\qquad \quad
\widehat{S}_{\widetilde{T}}(\widetilde{T}_j|\widehat{\bm{Z}}_j)^{-1}
-\widehat{S}_{\widetilde{T}}(\widetilde{T}_j|\bm{Z}_j)^{-1},\label{C.31}\\
&K_{m}\left(h_n^{-1}(\widehat{\bm{Z}}_j-\widehat{\bm{Z}})\right)- K_{m}\Big(h_n^{-1}(\bm{Z}_j-\bm{Z})\Big).\label{C.32}\end{align}

By means of reduction of fractions to a common denominator, the problem of finding the convergence rates of \eqref{C.31} translates into the rates of the following terms:  $$\widehat{f}(\widehat{\bm{Z}})-\widehat{f}(\bm{Z}) \quad,\quad
\widehat{S}_{\widetilde{T}}(\widetilde{T}_j|\widehat{\bm{Z}}_j)
-\widehat{S}_{\widetilde{T}}(\widetilde{T}_j|\bm{Z}_j).$$
Note that
$\widehat{f}(\widehat{\bm{Z}})-\widehat{f}(\bm{Z})=
\dfrac{1}{n}\sum\limits_{i=1}^{n}h_n^{-m}\left[K_{m}\left(h_n^{-1}(\bm{Z}_i-\widehat{\bm{Z}})\right)-K_{m}\Big(h_n^{-1}(\bm{Z}_i-\bm{Z})\Big)\right]$.
By the Taylor's expansion and assumption \eqref{A7}-\eqref{A8}, we have
\begin{align*}
&K_{m}\left(h_n^{-1}(\bm{Z}_i-\widehat{\bm{Z}})\right)-K_{m}\Big(h_n^{-1}(\bm{Z}_i-\bm{Z})\Big)\\
=&\nabla_{K_{m}}\big(h_n^{-1}(\bm{Z}_i-\bm{Z})\big)^Th_n^{-1}\left(\widehat{\bm{Z}}-\bm{Z}\right)+o_p\left(h_n^{-1}\|\widehat{\bm{Z}}-\bm{Z}\|\right)\\
=&O_p(n^{-b}h_n^{-1}).
\end{align*}
Thus,
\be\widehat{f}(\widehat{\bm{Z}})-\widehat{f}(\bm{Z})=O_p\left(n^{-b}h_n^{-(m+1)}\right).\label{C.33}\ee
We now evaluate $\widehat{S}_{\widetilde{T}}\big(\widetilde{T}_i|\widehat{\bm{Z}}_i\big)-\widehat{S}_{\widetilde{T}}\big(\widetilde{T}_i|\bm{Z}_i\big)$. We have
\begin{align}
&\widehat{S}_{\widetilde{T}}\big(\widetilde{T}_i|\widehat{\bm{Z}}_i\big)-\widehat{S}_{\widetilde{T}}\big(\widetilde{T}_i|\bm{Z}_i\big)\no\\
=&\dfrac{\dfrac{1}{n}\sum\limits_{j\ :\ \widetilde{T}_j>\widetilde{T}_i}h_n^{-m}K_{m}\big(h_n^{-1}(\bm{Z}_j-\widehat{\bm{Z}}_i)\big)}{\widehat{f}(\widehat{\bm{Z}}_i)}
-\dfrac{\dfrac{1}{n}\sum\limits_{j\ :\ \widetilde{T}_j>\widetilde{T}_i}h_n^{-m}K_{m}\big(h_n^{-1}(\bm{Z}_j-\bm{Z}_i)\big)}{\widehat{f}(\bm{Z}_i)}\no\\
=&\dfrac{\widehat{f}(\bm{Z}_i)n^{-1}\sum\limits_{j\ :\ \widetilde{T}_j>\widetilde{T}_i}h_n^{-m}K_{m}\big(h_n^{-1}(\bm{Z}_j-\widehat{\bm{Z}}_i)\big)-\widehat{f}(\widehat{\bm{Z}}_i)n^{-1}\sum\limits_{j\ :\ \widetilde{T}_j>\widetilde{T}_i}h_n^{-m}K_{m}\big(h_n^{-1}(\bm{Z}_j-\bm{Z}_i)\big)}
{\widehat{f}(\widehat{\bm{Z}}_i)\widehat{f}(\bm{Z}_i)}\no\\
=&\dfrac{\widehat{f}(\bm{Z}_i)n^{-1}\sum\limits_{j\ :\ \widetilde{T}_j>\widetilde{T}_i}h_n^{-m}\left[K_{m}\big(h_n^{-1}(\bm{Z}_j-\widehat{\bm{Z}}_i)\big)
-K_{m}\big(h_n^{-1}(\bm{Z}_j-\bm{Z}_i)\big)\right]+O_p\left(n^{-b}h_n^{-(m+1)}\right)}{\widehat{f}(\widehat{\bm{Z}}_i)\widehat{f}(\bm{Z}_i)}\no\\
=&O_p\left(n^{-b}h_n^{-(m+1)}\right)\label{C.34}.
\end{align}
as \eqref{C.32} satisfies that
\begin{align}&K_{m}\left(h_n^{-1}(\widehat{\bm{Z}}_j-\widehat{\bm{Z}})\right)-K_{m}\Big(h_n^{-1}(\bm{Z}_j-\bm{Z})\Big)\no\\
=&\nabla_{K_{m}}\big(h_n^{-1}(\bm{Z}_j-\bm{Z})\big)^Th_n^{-1}
\left[\widehat{\bm{Z}}_j-\bm{Z}_j-\left(\widehat{\bm{Z}}-\bm{Z}\right)\right]+
o_p\left(h_n^{-1}\left\|\widehat{\bm{Z}}_j-\bm{Z}_j-\left(\widehat{\bm{Z}}-\bm{Z}\right)\right\|\right)\no\\
=&O_p(n^{-b}h_n^{-1}).\label{zzz}
\end{align}
By \eqref{C.33}-\eqref{C.34}, we can obtain the convergence rates of  both terms in \eqref{C.31}:
\be\widehat{f}(\widehat{\bm{Z}})^{-1}-\widehat{f}(\bm{Z})^{-1}=-\dfrac{\widehat{f}(\widehat{\bm{Z}})-\widehat{f}(\bm{Z})}
{\widehat{f}(\widehat{\bm{Z}})\widehat{f}(\bm{Z})}=O_p\left(n^{-b}h_n^{-(m+1)}\right) \label{C.35}\ee
and
\be\widehat{S}_{\widetilde{T}}\big(\widetilde{T}_i|\widehat{\bm{Z}}_i\big)^{-1}-
\widehat{S}_{\widetilde{T}}\big(\widetilde{T}_i|\bm{Z}_i\big)^{-1}
=-\dfrac{\widehat{S}_{\widetilde{T}}\big(\widetilde{T}_i|\widehat{\bm{Z}}_i\big)-\widehat{S}_{\widetilde{T}}\big(\widetilde{T}_i|\bm{Z}_i\big)}
{\widehat{S}_{\widetilde{T}}\big(\widetilde{T}_i|\widehat{\bm{Z}}_i\big)\widehat{S}_{\widetilde{T}}\big(\widetilde{T}_i|\bm{Z}_i\big)}
=O_p\left(n^{-b}h_n^{-(m+1)}\right).\label{zzz1}\ee
\noindent {By \eqref{C.30}, \eqref{zzz}, and \eqref{zzz1}, it can be shown that
$$\text{\Rmnum{3}}_j=O_p\left(n^{-b}h_n^{-(2m+1)}\right),
$$
which, jointly with the law of large numbers, implies that
\be\text{\Rmnum{2}}_1-\text{\Rmnum{2}}_2=O_p\left(n^{-b}h_n^{-(2m+1)}\right).\label{C.36}\ee}
The assumption \eqref{A6} ensures the rationality of the convergence rate $O_p(n^{-b}h_n^{-(2m+1)})$.

Combining \eqref{C.27} and \eqref{C.28}, we have
\begin{align}&\dfrac{1}{n}\sum\limits_{i=1}^nI_i\left[\widehat{w}(\widetilde{T}_i,t,\bm{Z}_i)
-\widehat{w}(\widetilde{T}_i,t,\widehat{\bm{Z}}_i)\right]\no\\
=&\dfrac{1}{n}\sum\limits_{i=1}^nI_ie^{-\widehat{\Lambda}
(\widetilde{T}_i,t|\bm{Z}_i)}\left[\widehat{\Lambda}(\widetilde{T}_i,t|\widehat{\bm{Z}}_i)
-\widehat{\Lambda}(\widetilde{T}_i,t|\bm{Z}_i)\right]\cdot\Big(1+o_p(1)\Big)\no\\
=&\dfrac{1}{n}\sum\limits_{i=1}^nI_ie^{-\widehat{\Lambda}
(\widetilde{T}_i,t|\bm{Z}_i)}\left[\dfrac{\text{\Rmnum{2}}_1-\text{\Rmnum{2}}_2}{f(\widehat{\bm{Z}}_i)}
+\text{\Rmnum{2}}_1\left(\dfrac{1}{\widehat{f}(\widehat{\bm{Z}}_i)}
-\dfrac{1}{\widehat{f}(\bm{Z}_i)}\right)\right]\cdot\Big(1+o_p(1)\Big),\label{rrr1}
\end{align}
$\hfill{} \square$

%Denote
%$$\widehat{\bG}=\dfrac{1}{n}\Phi\widehat{\bw}_l=\dfrac{1}{n}\sum_{i=1}^{n}\bphi(\bm{X}_i)I(\widetilde{T}_i\geq t)+\dfrac{1}{n}\sum\limits_{i:\widetilde{T}_i<t,\triangle_i=0}\bphi(\bm{X}_i)\widehat{w}\big(\widetilde{T}_i,t,\widehat{\bm{Z}}_i\big).$$

Denote $$\widetilde{G}_1=\dfrac{1}{n}\sum_{i=1}^{n}\phi_1(\bm{X}_i)I(\widetilde{T}_i\geq t)+\dfrac{1}{n}\sum\limits_{i:\widetilde{T}_i<t,\triangle_i=0}\phi_1(\bm{X}_i)\widehat{w}\big(\widetilde{T}_i,t,\bm{Z}_i\big),$$

\begin{proposition}\label{prop3} Under {Assumptions \eqref{A1}-\eqref{A5}}, we have
\be\left|\widetilde{G}_1-E\big[\phi_1(\bm{X})I(T\geq t)\big]\right|=O_p(n^{-1/2}).\label{C.39}\ee
\end{proposition}

{\bf Proof}. In fact, under Assumptions \eqref{A1}-\eqref{A5}, its proof is very similar to the proof of (L1). By \eqref{C.14} and \eqref{C.19} we have
\begin{align*}&\dfrac{1}{n}\sum\limits_iI_i\phi_1(\bm{X}_i)\widehat{w}(\widetilde{T}_i,t,\bm{Z}_i)
=\dfrac{1}{n}\sum\limits_iI_i\phi_1(\bm{X}_i) e^{-\widehat{\Lambda}(\widetilde{T}_i,t|\bm{Z}_i)}\\
=&\dfrac{1}{n}\sum\limits_iI_i\phi_1(\bm{X}_i)\left[e^{-\Lambda_i}+e^{-\Lambda_i}
\left(f_i^{-1}n^{-1}\sum\limits_j\varepsilon_{ij}-{\tt Q}_{22}+{\tt Q}_{23}-{\tt Q}_{24}+O_p(n^{-1/2})\right)\right]\\
=&\dfrac{1}{n}\sum\limits_iI_i\phi_1(\bm{X}_i)w_i+n^{-2}\sum\limits_{i,j}I_i\phi_1(\bm{X}_i)w_if_i^{-1}\varepsilon_{ij}-n^{-1}\sum\limits_iI_i\phi_1(\bm{X}_i)w_i
\Big[{\tt Q}_{22}-{\tt Q}_{23}+{\tt Q}_{24}+O_p(n^{-1/2})\Big]\\
=&D_{31}+D_{32}-D_{33}.
\end{align*}
It is obvious that
$$D_{31}-ED_{31}=\dfrac{1}{n}\sum\limits_i\Big[I_i\phi_1(\bm{X}_i)w_i-E(I_i\phi_1\big(\bm{X}_i)w_i\big)\Big]=O_p(n^{-1/2}).$$
Similar to \eqref{C.21}, we can express $E(I_i\phi_1(\bm{X}_i)w_i)$ as
\begin{align}
  E(I_i\phi_1(\bm{X}_i)w_i) & = E\left[\phi_1(\bm{X})I(\widetilde{T}<t,\triangle=0)P(T\geq t|C,T\geq C,\bm{Z})\right]\no\\
   & =E\left[\phi_1(\bm{X})I(\widetilde{T}<t,\triangle=0)P(T\geq t|\widetilde{T},\triangle=0,\bm{Z})\right]\no\\
   & =E\left[\phi_1(\bm{X})I(\widetilde{T}<t,\triangle=0)E\left(I(T\geq t)\Big|\widetilde{T},\triangle=0,\bm{Z}\right)\right]\no\\
   & =E\left\{\phi_1(\bm{X})E\left[I(\widetilde{T}<t,\triangle=0)I(T\geq t)\Big|\widetilde{T},\triangle=0,\bm{Z}\right]\right\}\no\\
   & =E\left[\phi_1(\bm{X})I(\widetilde{T}<t,\triangle=0)I(T\geq t)\right]\no\\
   & =E\left[\phi_1(\bm{X})I(\widetilde{T}<t,\triangle=0,T\geq t)\right]\no\\
   & =E\Big[\phi_1(\bm{X})I(T\geq t,C<t)\Big]. \label{C.40}
\end{align}
The fifth ``$=$" of \eqref{C.40} is due to the conditional independence
$$(T,C)\independent \bm{X}\mid \bm{Z},$$
i.e., Condition 1.
The convergence rates of $D_{32}$ and $D_{33}$ can be obtained in
the same way as that of ${\tt Q}_{32}$ and ${\tt Q}_{33}$. Similarly, we obtain that
\be\dfrac{1}{n}\sum_iI_i\phi_1(\bm{X}_i)\widehat{w}(\widetilde{T}_i,t,\bm{Z}_i)=E\Big[\phi_1(\bm{X})I(T\geq t,C<t)\Big]+O_p(n^{-1/2}).\label{zzz2}\ee
By the central limit theorem, we can show that
$$\dfrac{1}{n}\sum\limits_{i=1}^{n}\phi_1(\bm{X}_i)I(\widetilde{T}_i\geq t)  =E\Big[\phi_1(\bm{X})I(T\geq t,C\geq t)\Big]+O_p(n^{-1/2}), $$
which, jointly with \eqref{zzz2}, yields \eqref{C.39}.

$\hfill{} \square$

{\bf Proof of Lemma \ref{lemma4}}

Note that $$\widehat{p}_l=\widehat{P}(T\geq t_l)-\widehat{P}(T\geq t_{l+1})=\dfrac{1}{n}\mathbf{1}_n^T\widehat{\bm{w}}_{\ell}.$$
To prove $(L1)$, we need only to show
$$\Big|\widehat{P}(T\geq t)-P(T\geq t)\Big|=O_p(n^{-\kappa}),$$
where
$$\widehat{P}(T\geq t)=\dfrac{1}{n}\sum\limits_{i=1}^{n}I(\widetilde{T}_i\geq t)+\dfrac{1}{n}\sum\limits_{i:\widetilde{T}_i<t,\triangle_i=0}\widehat{w}\big(\widetilde{T}_i,t,\widehat{\bm{Z}}_i\big).$$
By Proposition \ref{prop1}, it follows that $\Big|\widetilde{P}(T\geq t)-P(T\geq t)\Big|=O_p(n^{-1/2}).$, which implies that
\be\label{lll}|\widetilde{p}_l-p_l|=O_p(n^{-1/2}),\ee
where $\widetilde{p}_l=\widetilde{P}(T\geq t_l)-\widetilde{P}(T\geq t_{l+1})$.
In view of Proposition \ref{prop2}, we have
$$\Big|\widehat{P}(T\geq t)-\widetilde{P}(T\geq t)\Big|=O_p(n^{-\kappa})\ ,
\qquad 0<\kappa\leq \dfrac{1}{4},$$
which, jointly with \eqref{lll}, yields that
$$\left|{\widehat{p}_\ell}-p_\ell\right|= O_p(n^{-\kappa}).$$

Now we prove (L2). Denote
$$\widehat{\bG}=\dfrac{1}{n}\Phi\widehat{\bw}_l=\dfrac{1}{n}\sum_{i=1}^{n}\bphi(\bm{X}_i)I(\widetilde{T}_i\geq t)+\dfrac{1}{n}\sum\limits_{i:\widetilde{T}_i<t,\triangle_i=0}\bphi(\bm{X}_i)\widehat{w}\big(\widetilde{T}_i,t,\widehat{\bm{Z}}_i\big).$$
Similar to the proof of (L1), it suffices to show
\be\left\|\widehat{\bG}-E\left[\bphi(\bm{X})I(T\geq t)\right]\right\|=O_p(n^{-\kappa}).\label{fff}\ee
Recall that $\bphi(\bm{X})$ is a infinite dimensional vector
\be\bphi(\bm{X})=(\sqrt{a_1}\phi_1(\bm{X}),\sqrt{a_2}\phi_2(\bm{X}),\cdots)^T\in l_2.\label{C.37}\ee
Thus we evaluate every element of $\widehat{\bG}-E\left[\bphi(\bm{X})I(T\geq t)\right]$. Without loss of generality, we consider the first component $\widehat{G}_1$ given by
$$\widehat{G}_1=\dfrac{1}{n}\sum_{i=1}^{n}\phi_1(\bm{X}_i)I(\widetilde{T}_i\geq t)+\dfrac{1}{n}\sum\limits_{i:\widetilde{T}_i<t,\triangle_i=0}\phi_1(\bm{X}_i)\widehat{w}\big(\widetilde{T}_i,t,\widehat{\bm{Z}}_i\big).$$
We first show that
\be\left|\widehat{G}_1-E\big[\phi_1(\bm{X})I(T\geq t)\big]\right|=O_p(n^{-b}h_n^{-(2m+1)}).\label{C.38}\ee
Similar to the derivation of \eqref{rrr1}, we can obtain that
\begin{align*}
  \widetilde{G}_1-\widehat{G}_1 & = \dfrac{1}{n}\sum\limits_{i=1}^nI_i\phi_1(\bm{X}_i)\widehat{w}\left[\widehat{w}(\widetilde{T}_i,t,\bm{Z}_i)
-\widehat{w}(\widetilde{T}_i,t,\widehat{\bm{Z}}_i)\right]\\
&=\dfrac{1}{n}\sum\limits_{i=1}^nI_i\phi_1(\bm{X}_i)e^{-\widehat{\Lambda}
(\widetilde{T}_i,t|\bm{Z}_i)}\left[\dfrac{\text{\Rmnum{2}}_1-\text{\Rmnum{2}}_2}{f(\widehat{\bm{Z}}_i)}
+\text{\Rmnum{2}}_1\left(\dfrac{1}{\widehat{f}(\widehat{\bm{Z}}_i)}
-\dfrac{1}{\widehat{f}(\bm{Z}_i)}\right)\right]\cdot\Big(1+o_p(1)\Big)\\
&=O_p(n^{-b}h_n^{-(2m+1)}),
\end{align*}
{which, jointly with Proposition \ref{prop3}, implies \eqref{C.38}, i.e.,
$$\left|\widehat{G}_1-E\big[\phi_1(\bm{X})I(T\geq t)\big]\right|=O_p(n^{-b}h_n^{-(2m+1)}).$$ We choose $h_n=O(n^{-1/2d})$ and denote $\kappa=b-(m+1)/d$. The assumption \eqref{A6} ensures that $\kappa>0$ and hence we have
$$n^{\kappa}\left|\widehat{G}_1-E\big[\phi_1(\bm{X})I(T\geq t)\big]\right|=O_p(h_n),$$
which enable us to find a uniform upper bound of $n^{\kappa}\left|\widehat{G}_i-E\big[\phi_i(\bm{X})I(T\geq t)\big]\right|$ in probability. We denote this bound as $C$. Combined with \eqref{C.37}, it follows that
$$n^{\kappa}\left\|\widehat{\bG}-E\left[\bphi(\bm{X})I(T\geq t)\right]\right\|\leq\left(\sum_{i=1}^{\infty}a_i\right)^{1/2}C,$$
{which, jointly with \eqref{C.37} and the fact that $\sum_{i}a_i<\infty$} , concludes \eqref{fff}, i.e.,
$$\left\|\widehat{\bG}-E\left[\bphi(\bm{X})I(T\geq t)\right]\right\|=O_p(n^{-\kappa}),$$
and hence (L2) holds. $\hfill{}\square$}

\noindent{\bf D. Proof of Theorem 1}

\renewcommand{\theequation}{D.\arabic{equation}}
\setcounter{equation}{0}
Without loss of generality, we assume that $\widehat{\bm{\mu}}=(1/n)\sum\limits_{i=1}^n\bphi(\bm{X}_i)={\bf 0}$, for otherwise we
can subtract $\widehat{\bm{\mu}}$ from $\bphi(\bm{X}_i)$. From \eqref{nnnn} we have
\be\widehat{\Gamma}=\sum_{\ell=1}^L\widehat{p}_\ell\widehat{\bm{\mu}}_\ell\otimes\widehat{\bm{\mu}}_\ell
= \sum_{\ell=1}^L\dfrac{1}{\widehat{p}_{\ell}}\dfrac{1}{n}\Phi\widehat{\bw}_{\ell}
\otimes\big(\dfrac{1}{n}\Phi\widehat{\bw}_{\ell}\big).\label{D.1}\ee
Denote
\begin{align*}
  n_{\ell} & =\#\{T_i\in D_{\ell}:\ i=1,\cdots,n\},\ \ell=1,2,\cdots,L,\\
  \widehat{\bve}_{\ell} & =\left\{\left(\widehat{e}_{\ell}^{(1)},\cdots,\widehat{e}_{\ell}^{(n)}\right)^T
:\ \widehat{e}_{\ell}^{(i)}=I(n_1+\cdots+n_{\ell-1}+1\leq i\leq n_1+\cdots+n_{\ell}),\ i=1,\cdots,n\right\},\\
\end{align*}

\noindent let $\widetilde{\Gamma}$ be conditional covariance operator of $\bphi(\bm{X})$ corresponding to the
data without censoring, that is
\be\widetilde{\Gamma}=\sum_{\ell=1}^{L}\dfrac{1}{n_{\ell}/n}\dfrac{1}{n}\Phi\widehat{e}_{\ell}\otimes
\big(\dfrac{1}{n}\Phi\widehat{e}_{\ell}\big).\label{D.2}\ee
From \cite{Yao2003}, we have
\be\big\|\widehat{\Sigma}-\Sigma\big\|_{HS}=O_p\big(1/\sqrt{n}\big)\quad,\quad \big\|\widetilde{\Gamma}-\Gamma
\big\|_{HS}=O_p\big(1/\sqrt{n}\big).\label{D.3}\ee
To simplify the notion, we denote
\be\Xi=\Sigma^{-1}\Gamma\quad ,\quad \widehat{\Xi}=(\widehat{\Sigma}^2+\tau I)^{-1}\widehat{\Sigma}\widehat{\Gamma}\quad,
\quad\widetilde{\Xi}=(\widehat{\Sigma}^2+\tau I)^{-1}\widehat{\Sigma}\widetilde{\Gamma}$$
$$\Xi_1=(\widehat{\Sigma}^2+\tau I)^{-1}\Sigma\Gamma\quad,\quad\Xi_2=(\Sigma^2+\tau I)^{-1}\Sigma\Gamma,\label{D.4}\ee
then
\be\|\widehat{\Xi}-\Xi\|_{HS}\leq \|\Pi_1\|_{HS}+\|\Pi_2\|_{HS}+\|\Pi_3\|_{HS}+\|\Pi_4\|_{HS},\label{D.5}\ee
where
\be\Pi_1=\widehat{\Xi}-\widetilde{\Xi}\quad,\quad\Pi_2=\widetilde{\Xi}-\Xi_1
\quad,\quad\Pi_3=\Xi_1-\Xi_2\quad,\quad\Pi_4=\Xi_2-\Xi.\label{D.6}\ee
Now we evaluate each of $\Pi_1\sim\Pi_4$. We fiirst consider $\Pi_1$, and we have
\be\|\Pi_1\|_{HS}\leq \big\|(\widehat{\Sigma}^2+\tau I)^{-1}\big\| \ \big\|\widehat{\Sigma}\big\|\ \big\|\widehat{\Gamma}-\widetilde{\Gamma}\big\|_{HS}.\label{D.7}\ee
The rank of operator $\widehat{\Sigma}$ is no more than $n-1$, so the smallest eigenvalue of $\widehat{\Sigma}$ is 0. Taking this into account, we have
\be\big\|(\widehat{\Sigma}^2+\tau I)^{-1}\big\|=\dfrac{1}{\tau},\label{D.eigen}\ee
which is the largest eigenvalue of $(\widehat{\Sigma}^2+\tau I)^{-1}$, in other words, the reciprocal
of the smallest eigenvalue of $(\widehat{\Sigma}^2+\tau I)$.
Note that
\begin{align}
  \widehat{\Gamma}-\widetilde{\Gamma} & =\sum_{\ell=1}^L\dfrac{1}{\widehat{p}_{\ell}}\dfrac{1}{n}\Phi\widehat{\bw}_{\ell}
\otimes\big(\dfrac{1}{n}\Phi\widehat{\bw}_{\ell}\big)\no\\
    & = \sum_{\ell=1}^{L}(\Pi_{\ell 1}+\Pi_{\ell 2}+\Pi_{\ell 3}),\label{D.8}
\end{align}
where
\begin{align}
  \Pi_{l1} & =\left(\dfrac{1}{\widehat{p}_{\ell}}-\dfrac{1}{n_{\ell}/n}\right)\dfrac{1}{n}\Phi\widehat{\bw}_{\ell}
\otimes\big(\dfrac{1}{n}\Phi\widehat{\bw}_{\ell}\big), \no\\
  \Pi_{l2} & =\dfrac{1}{n_{\ell}/n}\dfrac{1}{n}\Phi\widehat{\bw}_{\ell}\otimes \left(\dfrac{1}{n}\Phi\widehat{\bw}_{\ell}-\dfrac{1}{n}\Phi\widehat{e}_{\ell}\right), \label{D.9}\\
  \Pi_{l3} & =\dfrac{1}{n_{\ell}/n}\left(\dfrac{1}{n}\Phi\widehat{\bw}_{\ell}-\dfrac{1}{n}\Phi\widehat{e}_{\ell}\right)
  \otimes (\dfrac{1}{n}\Phi\widehat{e}_{\ell}).\no
\end{align}

\noindent By (L1) in Lemma \ref{lemma4}, it follows that
\be\|\Pi_{l1}\|_{HS}=\left|\dfrac{1}{\widehat{p}_{\ell}}-\dfrac{1}{n_{\ell}/n}\right|
\left\|\dfrac{1}{n}\Phi\widehat{\bw}_{\ell}\otimes\big(\dfrac{1}{n}\Phi\widehat{\bw}_{\ell}\big)\right\|_{HS}=O_p(n^{-\kappa}).\label{D.10}\ee
Note that $\big\|\frac{1}{n}\Phi\widehat{\bm{w}}_\ell-E[\bphi(\bX)I(T\in D_\ell)]\big\|=O_p(n^{-\kappa})$. Then by (L2) in Lemma \ref{lemma4}, we have
\begin{align}
    & \|\Pi_{l2}\|_{HS} \no\\
     \leq& \left|\dfrac{1}{n_{\ell}/n}\right|\Big\|\dfrac{1}{n}\Phi\widehat{\bw}_{\ell}\Big\|
     \left(\Big\|\frac{1}{n}\Phi\widehat{\bm{w}}_\ell-E\big[\bphi(\bX)I(T\in D_\ell)\big]\Big\|
     +\Big\|\frac{1}{n}\Phi\widehat{e}_\ell-E\big[\bphi(\bX)I(T\in D_\ell)\big]\Big\|\right)\no\\
   =&O_p(n^{-\kappa}),\label{D.11}
\end{align}
and
\begin{align}
    & \|\Pi_{l3}\|_{HS} \no\\
     \leq& \left|\dfrac{1}{n_{\ell}/n}\right|\Big\|\dfrac{1}{n}\Phi\widehat{e}_{\ell}\Big\|
     \left(\Big\|\frac{1}{n}\Phi\widehat{\bm{w}}_\ell-E\big[\bphi(\bX)I(T\in D_\ell)\big]\Big\|
     +\Big\|\frac{1}{n}\Phi\widehat{e}_\ell-E\big[\bphi(\bX)I(T\in D_\ell)\big]\Big\|\right)\no\\
    =&O_p(n^{-\kappa}).\label{D.12}
\end{align}
Combining \eqref{D.8}-\eqref{D.12}, we obtain that
\be\big\|\widehat{\Gamma}-\widetilde{\Gamma}\big\|_{HS}\leq \sum_{\ell=1}^{L}(\|\Pi_{l1}\|_{HS}+
\|\Pi_{l2}\|_{HS}+\|\Pi_{l3}\|_{HS})=O_p(n^{-\kappa}).\label{D.13}\ee
 {which, jointly with \eqref{D.eigen}, implies that \be\|\Pi_1\|_{HS}=O_p\Big(\dfrac{1}{\tau n^{\kappa}}\Big).\label{D.14}\ee}
Next we evaluate $\Pi_2$. {By \eqref{D.eigen}, it follows that
\be\|\Pi_2\|_{HS}\leq \big\|(\widehat{\Sigma}^2+\tau I)^{-1}\big\|\left(\big\|\widehat{\Sigma}\big\|
\big\|\widehat{\Gamma}-\Gamma\big\|_{HS}+\big\|\widehat{\Sigma}-\Sigma\big\|_{HS}\big\|\Gamma\big\|\right)
=O_p\Big(\dfrac{1}{\tau\sqrt{n}}\Big).\label{D.15}\ee}

To evaluate $\Pi_3$, we use the spectral decomposition of $\Gamma$ and obtain that
\be\Xi_1=\sum_{j=1}^{d_{\Gamma}}\zeta_j\Big(\big(\widehat{\Sigma}^2+\tau I\big)^{-1}\Sigma\varphi_j\Big)\otimes\varphi_j\quad,\quad
\Xi_2=\sum_{j=1}^{d_{\Gamma}}\zeta_j\Big(\big(\Sigma^2+\tau I\big)^{-1}\Sigma\varphi_j\Big)\otimes\varphi_j,\label{D.16}\ee
and
\be\left(\big(\widehat{\Sigma}^2+\tau I\big)^{-1}-\big(\Sigma^2+\tau I\big)^{-1}\right)\Sigma\varphi_j
=\big(\widehat{\Sigma}^2+\tau I\big)^{-1}\big(\Sigma^2-\widehat{\Sigma}^2\big)\big(\Sigma^2+\tau I\big)^{-1}\Sigma^2\widetilde{\varphi}_j,\label{D.17}\ee
which yields that
\be\|\Pi_3\|_{HS}\leq \sum_{j=1}^{d_{\Gamma}}\zeta_j\big\|(\widehat{\Sigma}^2+\tau I)^{-1}\big\|\
\big\|\widehat{\Sigma}^2-\Sigma^2\big\|_{HS}\ \big\|(\Sigma^2+\tau I)^{-1}\Sigma^2\big\|\ \|\widetilde{\varphi}_j\|
=O_p\Big(\dfrac{1}{\tau\sqrt{n}}\Big).\label{D.18}\ee

As for $\Pi_4$, it can be seen that
\be\|\Pi_4\|_{HS}=\sum_{j=1}^{d_{\Gamma}}\zeta_j\Big\|\left(\big(\Sigma^2+\tau I\big)^{-1}\Sigma-\Sigma^{-1}\right)\varphi_j\Big\|,\label{D.19}\ee
and for each $j=1,2,\cdots,d_{\Gamma}$,
\begin{align}
    & \Big\|\left(\big(\Sigma^2+\tau I\big)^{-1}\Sigma-\Sigma^{-1}\right)\varphi_j\Big\|\leq
    \Big\|\big(\Sigma^2+\tau I\big)^{-1}\Sigma^2\widetilde{\varphi}_j-\widetilde{\varphi}_j\Big\| \no\\
  = & \left\|\sum\limits_{i=1}^{\infty}\left(\dfrac{\nu_i^2}{\tau+\nu_i^2}-1\right)\langle \widetilde{\varphi}_{j},\psi_i\rangle \psi_i\right\|=
  \left(\sum\limits_{i=1}^{\infty}\dfrac{\tau^2}{(\tau+\nu_i^2)^2}\langle \widetilde{\varphi}_{j},\psi_i\rangle^2\right)^{1/2}\no\\
  \leq & \dfrac{\tau}{\nu_N^2}\left(\sum\limits_{i=1}^{N}\langle \widetilde{\varphi}_{j},\psi_i\rangle^2\right)^{1/2}
+\left(\sum\limits_{i=N+1}^{\infty}\langle \widetilde{\varphi}_{j},\psi_i\rangle^2\right)^{1/2}\no\\
  = & \dfrac{\tau}{\nu_N^2}\big\|\Psi_N(\widetilde{\varphi}_j)\big\|+\big\|\Psi_N^{\perp}(\widetilde{\varphi}_j)\big\|,\label{D.20}
\end{align}
{which, joint with \eqref{D.14}, \eqref{D.15}, \eqref{D.18}, and \eqref{D.19}, implies that \eqref{d3.5} holds.}

%\Acknowledgements{This work was supported by National Natural Science Foundation of China (Grant No. \bXX\bXX\bXX).}


\begin{thebibliography}{99}
\bibitem{Alfons2017}Alfons A, Croux C, Filzmoser P. Robust maximum association estimators. Journal of the American Statistical Association, 2017, 112(517): 436-445.
\bibitem{Alfons2016}Alfons A, Croux C, Filzmoser P. Robust maximum association between data sets: The R package ccaPP. Austrian Journal of Statistics, 2016, 45(1): 71-79.
\bibitem{RKHSfirst}Aronszajn N. Theory of reproducing kernels. T Am Math Soc, 1950: 337-404.
\bibitem{Berlinet2011}Berlinet A, Thomas-Agnan C. Reproducing kernel Hilbert spaces in probability and statistics. Springer Science \& Business Media, 2011.
\bibitem{Baker1981}Baker, C., \& McKeague, I. (1981). Compact Covariance Operators. Proceedings of the American Mathematical Society, 83(3), 590-593. doi:10.2307/2044126
\bibitem{Dantzig2007}Candes E, Tao T. The Dantzig selector: Statistical estimation when p is much larger than n. The annals of Statistics, 2007, 35(6): 2313-2351.
\bibitem{COOK1998}Cook R D. Regression Graphics: Ideas for Studying Regressions through Graphics. Wiley, New York, 1998, 73: 215.
\bibitem{COOK1999}Cook R D, Lee H. Dimension reduction in binary response regression. Journal of the American Statistical Association, 1999, 94(448): 1187-1200.
\bibitem{COOK1991}Cook R D, Weisberg S. Sliced inverse regression for dimension reduction: Comment. Journal of the American Statistical Association, 1991, 86(414): 328-332.
\bibitem{Cui15}Cui Wenquan, Wu Chenglong. An approach to estimating nonlinear sufficient dimension reduction subspace for censored survival data. Journal of University of Science and Technology of China, 2015,45(9):709-716.
\bibitem{Duan2003}Duan K, Keerthi S S, Poo A N. Evaluation of simple performance measures for tuning SVM hyperparameters. Neurocomputing, 2003, 51: 41-59.
\bibitem{spline}De Boor, C., De Boor, C., Mathématicien, E. U., De Boor, C., \& De Boor, C. (1978). A practical guide to splines (Vol. 27, p. 325). New York: springer-verlag.
\bibitem{Diaconis1984}Diaconis P, Freedman D. Asymptotics of graphical projection pursuit. The annals of statistics, 1984: 793-815.
\bibitem{scad2}Fan J, Li R. \text{var}iable selection via nonconcave penalized likelihood and its oracle properties. J Am Stat Assoc, 2001, 96: 1348-1360.
\bibitem{SIS1}Fan J, Lv J. Sure independence screening for ultrahigh dimensional feature space. J R Stat Soc B, 2008, 70: 849-911.
\bibitem{SIS2}Fan J, Song R. Sure independence screening in generalized linear models with NP-dimensionality. Ann Stat, 2010, 38: 3567-3604.
\bibitem{Yao2003}Ferr\'e, L., Yao A F. Functional sliced inverse regression analysis. Statistics, 2003, 37(6): 475-488.
\bibitem{Fukumizu09}Fukumizu K, Bach F R, Jordan M I. Kernel dimension reduction in regression. Ann Stat, 2009: 1871-1905.
\bibitem{multikernel2011}Gönen, M., \& Alpaydın, E. (2011). Multiple kernel learning algorithms. Journal of machine learning research, 12(Jul), 2211-2268.
\bibitem{Johnson02}Johnson R A, Wichern D W. Applied multivariate statistical analysis. Upper Saddle River, NJ: Prentice hall, 2002.
\bibitem{Jones1993}Jones M C, Foster P J. Generalized jackknifing and higher order kernels. Journal of Nonparametric Statistics, 1993, 3(1): 81-94.
\bibitem{Kato2013} Kato T. Perturbation theory for linear operators. Springer Science \& Business Media, 2013.
\bibitem{Keerthi2003}Keerthi, S. Sathiya, and Chih-Jen Lin. "Asymptotic behaviors of support vector machines with Gaussian kernel." Neural Comput 15.7 2003: 1667-1689.
\bibitem{GSIRGSAVE}Lee K Y, Li B, Chiaromonte F. A general theory for nonlinear sufficient dimension reduction: Formulation and estimation. The Annals of Statistics, 2013, 41(1): 221-249.
\bibitem{directionalregression}Li B, Wang S. On directional regression for dimension reduction. J Am Stat Assoc, 2007, 102: 997-1008.
\bibitem{functional2017}Li, B., \& Song, J. (2017). Nonlinear sufficient dimension reduction for functional data. The Annals of Statistics, 45(3), 1059-1095.
\bibitem{Li1991}Li K C. Sliced inverse regression for dimension reduction. J Am Stat Assoc, 1991, 86: 316-327.
\bibitem{Li1999}Li K C, Wang J L, Chen C H. Dimension reduction for censored regression data. Ann Stat, 1999, 27: 1-23.
\bibitem{LiLi04}Li L, Li H. Dimension reduction methods for microarrays with application to censored survival data. Bioinformatics, 2004, 20: 3406-3412.
\bibitem{Li2008}Li L, Yin X. Sliced inverse regression with regularizations. Biometrics, 2008, 64: 124-131.
\bibitem{lungdata}Loprinzi, C. L., Laurie, J. A., Wieand, H. S., Krook, J. E., Novotny, P. J., Kugler, J. W., ... \& Klatt, N. E. (1994). Prospective evaluation of prognostic variables from patient-completed questionnaires. North Central Cancer Treatment Group. Journal of Clinical Oncology, 12(3), 601-607.
\bibitem{lu2011}Lu W, Li L. Sufficient dimension reduction for censored regressions. Biometrics, 2011, 67: 513-523.
\bibitem{MaZhu12}Ma, Y. and Zhu, L. (2012). A semiparametric approach to dimension reduction. J. Amer. Statist. Assoc., 107(497),168-179.
\bibitem{MaZhu13}Ma Y, Zhu L. A review on dimension reduction. Int Stat Rev, 2013, 81: 134-150.
\bibitem{wavelets2012}Ogden, T. (2012). Essential wavelets for statistical applications and data analysis. Springer Science \& Business Media.
\bibitem{scaledlasso}Sun, T., \& Zhang, C. H. (2012). Scaled sparse linear regression. Biometrika, 99(4), 879-898.
\bibitem{Bernhard02}Scholkopf B, Smola A J. Learning with kernels: support vector machines, regularization, optimization, and beyond. MIT press, 2001.
\bibitem{MayaStephan14}Shevlyakova M, Morgenthaler S. Sliced inverse regression for survival data. Stat Pap, 2014, 55: 209-220.
\bibitem{multikernel2006}Sonnenburg, S., Rätsch, G., Schäfer, C., \& Schölkopf, B. (2006). Large scale multiple kernel learning. Journal of Machine Learning Research, 7(Jul), 1531-1565.
\bibitem{LASSO}Tibshirani R. Regression shrinkage and selection via the lasso. J R Stat Soc B, 1996: 267-288.
\bibitem{wavelets2009}Vidakovic, B. (2009). Statistical modeling by wavelets (Vol. 503). John Wiley \& Sons.
\bibitem{Xuerong10}Wen X M. On sufficient dimension reduction for proportional censorship model with covariates. Comput Stat Data An, 2010, 54: 1975-1982.
\bibitem{Wahba1990}Wahba G. Spline models for observational data. Siam, 1990.
\bibitem{Wu2008}Wu H M. Kernel sliced inverse regression with applications to classification. J Comput Graph Stat, 2008, 17: 590-610.
\bibitem{Wuqiang2013}Wu Q, Liang F, Mukherjee S. Kernel sliced inverse regression: Regularization and consistency.  Abstr Appl Anal. 2013, Special Issue 2013: 1-11.
\bibitem{MAVE}Xia Y, Tong H, Li W K, et al. An adaptive estimation of dimension reduction space. J R Stat Soc B, 2002, 64: 363-410.
\bibitem{dMAVE}Xia Y. A constructive approach to the estimation of dimension reduction directions. The Annals of Statistics, 2007, 35(6): 2654-2690.
\bibitem{XIA2011}Xia Y, Zhang D, Xu J. Dimension reduction and semiparametric estimation of survival models. J Am Stat Assoc, 2010, 105: 278-290.
\bibitem{Yeh09}Yeh Y R, Huang S Y, Lee Y J. Nonlinear dimension reduction with kernel sliced inverse regression. IEEE T Knowl Data En, 2009, 21: 1590-1603.
\bibitem{grouplasso}Yuan, M., \& Lin, Y. (2006). Model selection and estimation in regression with grouped variables. Journal of the Royal Statistical Society: Series B (Statistical Methodology), 68(1), 49-67.
\bibitem{RSIR1991}Zhong W, Zeng P, Ma P, et al. RSIR: regularized sliced inverse regression for motif discovery. Bioinformatics, 2005, 21: 4169-4175.
\bibitem{kernelIR}Zhu L X, Fang K T. Asymptotics for kernel estimate of sliced inverse regression. The Annals of Statistics, 1996, 24(3): 1053-1068.
\bibitem{ALASSO} Zou, H. (2006). The Adaptive lasso and Its Oracle Properties. \emph{Journal of the American Statistical Association}, \textbf{101}, 1418-1429.

\bibitem{LASSOnet}Zou H, Hastie T. Regularization and variable selection via the elastic net. J R Stat Soc B, 2005, 67: 301-320.

\end{thebibliography}
\end{document}